\documentclass[11pt]{article}
\usepackage[dvips]{graphics}
\usepackage{amsmath}
\usepackage{amssymb}
\usepackage{amscd}
\usepackage{amsthm}
\usepackage{amsopn}
\usepackage{xspace}
\usepackage{verbatim}
\usepackage{amsmath}
\usepackage{amsfonts}
\usepackage{colortbl}
\usepackage{hyperref}
\usepackage{amsmath}
\usepackage{color}
\usepackage{srcltx}
\usepackage[utf8]{inputenc}
\usepackage{blindtext}
\usepackage[a4paper, total={6in, 10in}]{geometry}

\newcommand{\C}{\mathbb{C}}
\newcommand{\N}{\mathbb{N}}
\newcommand{\R}{\mathbb{R}}

\newcommand\AMSname{AMS subject classifications}

\newtheorem{prop}{Proposition}[section]
\newtheorem{thm}[prop]{Theorem}

\newtheorem{lemme}[prop]{Lemma}
\newtheorem{remarque}[prop]{Remark}

\newenvironment{preuve}[1][]{\noindent {\it Proof #1}  : } {\hbox{~}\qed}

\title{

Existence results for problems involving nonlocal operator with an asymmetric weight and with  a critical nonlinearity  }
\author{Sana Benhafsia \thanks{Université Paris-Est Créteil, LAMA, Laboratoire d'Analyse et
de Mathématiques Appliquées, CNRS UMR 8050, UPEC, F-94010 Créteil, France. E-mail:
sana.ben-hafsia@u-pec.fr} and $\,$  Rejeb Hadiji\thanks{Université Paris-Est Créteil, LAMA, Laboratoire d'Analyse et
de Mathématiques Appliquées, CNRS UMR 8050, UPEC, F-94010 Créteil, France. E-mail:
rejeb.hadiji@u-pec.fr  }}
\begin{document}
\date{}

\maketitle

\begin{abstract}
Recently, a great attention has been focused on the study of fractional and nonlocal operators of elliptic type, both for
the pure mathematical research and in view of concrete real-world applications.
We consider the following nonlocal problem on the  $\mathbb{H}_0^s(\Omega) \subset L^{q_s}(\Omega)$,  with $q_s :=\frac{2n}{n-2s}$,  $s\in ]0, 1[$ and  $n\geq 3$
\begin{equation*}\label{problème général}
\inf_{\substack{u\in \mathbb{H}_0^s(\Omega) \\ ||u||_{L^{q_s}(\mathbb{R}^n)}=1}} \left\{ \displaystyle\int_{\mathbb{R}^n}p(x) \bigg(\displaystyle\int_{\mathbb{R}^n}\frac{|u(x)-u(y)|^2}{|x-y|^{n+2s}}dy\bigg)dx-\lambda \displaystyle\int_\Omega |u(x)|^q dx\right\},
\end{equation*}
where $\Omega$ is a bounded domain in $\mathbb{R}^n,   p :\mathbb{R}^n \rightarrow \R$ is a given positive weight presenting a global positive
minimum $p_0 >0$ at
$a \in \Omega$ and $\lambda$ is a real constant. After introducing the fractional laplacian $(-\Delta)^s$, we will  show that for $q=2$ the infimum above is achieved on  the set
$\{u\in \mathbb{H}_0^s(\Omega), ||u||_{L^{q_s}(\Omega)}=1\}$  for a class of $k,  s,  \lambda$ and $n$ and for $q\geq 2$, the problem has also  a non-ground state solution.
\end{abstract}


\noindent2010 \AMSname: 35J20, 35J25, 35H30, 35J60.

\section{Introduction,  notations and statement of the results}
\subsection{The fractional nonlinear problem with weight and its relation with the ordinary nonlinear problem}
Let $\Omega$ be a smooth bounded domain of $\mathbb{R}^n$,  $n \geq 3$.
We are interested in the following nonlinear problem involving the fractional Laplacian, for $u \in \mathbb{H}_0^s(\Omega)$
\begin{equation}\label{formulation faible}
\begin{array}{ll}
\displaystyle\int_{\mathbb{R}^n\times \mathbb{R}^n}p(x)\frac{(u(x)-u(y))(\varphi(x)-\varphi(y))}{|x-y|^{n+2s}}dxdy -\lambda \displaystyle\int_{\Omega} |u(x)|^{q-2} u(x) \varphi(x) dx \\  = \displaystyle\int_{\Omega} |u(x)|^{q_s-2}u(x)\varphi(x) dx
\end{array}
\end{equation}
for any $\varphi \in \mathbb{H}_0^s(\Omega)$, where the space $\mathbb{H}_0^s(\Omega)$ is defined by
\begin{equation}
\mathbb{H}_0^s(\Omega) :=\big\{u\in L^2(\Omega), \frac{|u(x)-u(y)|}{|x-y|^{\frac{n}{2}+s}} \in L^2(\mathbb{R}^n \times \mathbb{R}^n \setminus \Omega^c \times \Omega^c), u(x)= 0,\forall x\in \mathbb{R}^n \setminus{\Omega }\big\}
\end{equation}
where $2\leq q<q_s$, $q_s$ is the critical fractional Sobolev exponent $q_s :=\frac{2n}{n-2s}$, $\lambda >0$. Note that it is well known that the embedding $\mathbb{H}_0^s(\Omega) \hookrightarrow L^r(\Omega)$ is continuous for any $1\leq r \leq\frac{2n}{n-2s}$. Moreover this embedding is compact for $1 \leq r <\frac{2n}{n-2s}$, see (\cite{BRS}, lemma 1.31). We assume that $p : \mathbb{R}^n \rightarrow \R$ 
 is a positive weight in $C(\mathbb{R}^n)$, we assume also  that 
 $p$ represents a global positive minimum $p_0$ at $a\in \Omega$ and satisfies in  $ B(a,4\eta)$, $\eta >0$, $\kappa >0$ and $k>1$
\begin{equation}\label{def de p}
p(x)\leq p_0 + \kappa |x-a|^k,
\end{equation}
and
 \begin{equation}\label{p-p_0 dans L1}
 p-p_0 \in L^1(\mathbb{R}^n).
 \end{equation}
We remind that for $p=1$, \eqref{formulation faible} is the weak formulation  of the following problem
\begin{equation}
      \left\{
\begin{aligned}
\begin{array}{ll}
(-\Delta)^s u - \lambda u^{q-2} u= |u|^{q_s-2} u&  \hbox{in }   \Omega \\
u=0 &  \hbox{in} \ \mathbb{R}^n\setminus\Omega, 
\end{array}\end{aligned} \right.
\end{equation}
which was studied in \cite{BRS} when $q=2, n\geq 4s$ and $\lambda \in ]0, \lambda_{1, s}[ $ where $\lambda_{1, s}$ denotes the first eigenvalue of the nonlocal operator $(-\Delta)^s$ with homogeneous Dirichlet boundary datum.
In this paper, we will break the symmetry of the problem by introducing the weight $p$ as described before. For the first main result we didn't require that $p$ is bounded, but for the second one, $p$ should be bounded. We extend the result which was already done in the ordinary Laplacian in \cite{HY} to the case of the fractional Laplacian. Other authors gave a basic introduction to the fractional Laplacian operator, see \cite{NPV} and references therein, see also \cite{AB}, \cite{ABS}, \cite{BRS} and \cite{NS} where authors dealt with non local fractional problems. Some other authors studied variants of the fractional Laplacian, see \cite{GQ}, \cite{LP} and \cite{CS}.\\
   Let us take $p$ defined as previously in \eqref{def de p},  $s\in]0, 1[$ and $n\geq 3$, we define the infimum $S_{s,\lambda}(p)$ by
  \begin{equation}\label{infimum}
S_{s, \lambda}(p):=\inf_{\substack{u\in \mathbb{H}_0^s(\Omega) \\ ||u||_{L^{q_s}(\mathbb{R}^n)}=1}} \left\{ \displaystyle\int_{\mathbb{R}^n}p(x) \bigg(\displaystyle\int_{\mathbb{R}^n}\frac{|u(x)-u(y)|^2}{|x-y|^{n+2s}}dy\bigg)dx-\lambda \displaystyle\int_\Omega |u(x)|^2dx \right\} .
\end{equation}
 The study of the  infimum \eqref{infimum},  shows that the existence of minimizers depends, not only on the behavior of 
 $p$ near its global minimum but also on its behavior at $+\infty$.

The well known fractional Sobolev inequalities were first considered in a remarkable paper by Lieb in  \cite{L}; see also \cite{FLS}; \cite{CT} or the survey \cite{NPV}.

 Let us compare this  work to what is known  in the literature concerning problems related to the Yamabe problem in the local case
\begin{equation}
S_\lambda(p)=\inf_{\substack{u\in H_0^1(\Omega)\\ ||u||_{L^{q}(\Omega)}=1}}\displaystyle\int_{\Omega}p(x)|\nabla u(x)|^2dx-\lambda \displaystyle\int_\Omega |u(x)|^2dx, 
\end{equation}
where $q=\frac{2n}{n-2}$ is the critical exponent for the Sobolev embedding $H_0^1(\Omega) \subset L^q(\Omega)$. 
In \cite{HY},  the authors treated the case where the problem presents a positive weight with the ordinary Laplacian ($s=1$). They proved in particular,  the existence of minimizers of $S_\lambda(p)$ for $0<\lambda<\lambda_1^{div}$ if $n\geq4$ and $k>2$,  and for 
$\lambda^*<\lambda<\lambda_1^{div}$ if $n\geq3$ and $0<k<2$, $k=2$ is critical for the problem,   and in other subcases which are well detailed in (\cite{HY},Theorem 1.1) with $k$ is a positive constant that appears in the expression of the weight $p$,  $\lambda_1^{div}$ is the first eigenvalue of $-div(p(x)\nabla.)$ on $\Omega$ with zero Dirichlet boundary condition and $\lambda^*$ is a positive constant. The method used for the proof of this result is,  first to show that $S_\lambda(p) <p_0 S$,  with $S$ is the best Sobolev constant defined by
\begin{equation}
S :=\inf_{\substack{u\in H_0^1(\Omega)\\||u||_{L^q(\Omega)}=1}}\displaystyle\int_{\Omega}|\nabla u(x)|^2dx,
\end{equation}
then,  they prove that the infimum $S_\lambda(p)$ is achieved. 
In the same context, similar questions were studied in \cite{H}  where the author investigate the problem with a weight and a nonvanishing boundary datum and in   \cite{HV} where authors dealt with  a nonlinear eigenvalue problem with a variable weight.

There are extensive works that study properties of fractional Laplacian operators as non-local operators together with applications to free-boundary value problems and non-local minimal surfaces by Caffarelli and Silvestre and many
others.\\
 Frank, Lieb and Seirenger treated in \cite{FLS} weighted cases namely $\displaystyle\frac{dx dy}{|x|^\alpha |y|^\beta}$ in order to prove Hardy-Lieb-Thirring inequalities. \\
Chang and González \cite{CG}, formulated a fractional $s$-Yamabe problem that include the boundary Yamabe problem studied by Escobar. They highlight a Hopf-type maximum principle together with interplay between analysis of weighted trace Sobolev inequalities and conformal structure of the underlying manifolds and they obtained some properties for the fractional case that are analogous to the original Yamabe problem.

The problem involving the ordinary Laplacian in the case where $p$ is a constant was originally studied by Aubin in  \cite{A}, Brezis-Nirenberg in  \cite{BN1}, \cite{B}. 

In \cite{BH}, we studied  polyharmonic operator, namely $(-\Delta)^r$ where $r\in \mathbb{N} $, $ r\geq 2$ with critical Sobolev exponent. We refer to  \cite{LP} for a complete history of the problem and geometrical motivations.\\
Servadei and Valdinocci dealt in \cite{SV} with a non-ground state solutions in the fractional Laplacian case without weight. The solutions are constructed with a variational method by a min-max procedure on the associated energy functional.

\noindent I refer to \cite{HY} and \cite{HMPY} for the study of existence and multiplicity solutions of above problem in the case where $s =1$   with the presence of a  smooth and positive  weight and  with homogeneous Dirichlet boundary condition. Nevertheless, in \cite{HY}, it is shown that if $p$ has  discontinuities then  a solution of $S_0 (p,0)$ still exists.

\noindent For more general weights, depending on $x$ and on $u$, in two papers written with Bae, Yazidi \cite{BHY} and  with Vigneron \cite{HV}, we  showed  that in the case of homogeneous Dirichlet  boundary condition and in the presence of a linear perturbation the corresponding minimizing problem possesses a solution. The model of the  weight is $p(x,u) = \alpha + \vert x\vert^\beta \vert u\vert^k$ with positive parameters $\alpha$, $\beta$ and $k$. Note that in this case natural scalings appear and the answer depends on the ratio ${\beta \over k}$. I refer to  \cite{LP} for a complete history of the problem and geometrical motivations.

In \cite{GQ} the authors formulated a fractional $s$-Yamabe problem that include the boundary Yamabe problem studied by Escobar, see \cite{Es}. They highlight a Hopf-type maximum principle together with interplay between analysis of weighted trace Sobolev inequalities and conformal structure of the underlying manifolds and they obtained some properties for the fractional case that are analogous to the original Yamabe problem, see \cite{Es}.
\subsection{Some definitions} 
One of the aims of this paper is to study non local problems driven by $(-\Delta)^s$ (or its generalization) and with Dirichlet boundary data via variational methods. For this purpose,  we need to work in a suitable fractional Sobolev space: for this,  we consider a functional analytical setting that is inspired by (but not equivalent to) the fractional Sobolev spaces in order to correctly encode the Dirichlet boundary datum in the
variational formulation. This section is devoted to the definition of this space as well
as to its properties. Therefore, before setting the main result,  we start by defining the fractional Laplacian and the fractional Sobolev spaces.\\
Let $s\in ]0, 1[$ and let $S(\mathbb{R}^n)$ be the Schwartz space of rapidly decaying $C^\infty$ functions in $\mathbb{R}^n$. We define the non local operator $(-\Delta)^s: S(\mathbb{R}^n) \longrightarrow  L^2(\mathbb{R}^n)$ by 
\begin{equation}
   \begin{array}{cl}
(-\Delta)^s u(x)   & := c(n, s) P.V.\displaystyle\int_{\mathbb{R}^n} \frac{u(x)-u(y)}{|x-y|^{n+2s}}dy \\
         & = c(n, s) \displaystyle\lim_{\varepsilon \rightarrow{0^+} }\displaystyle\int_{\mathbb{R}^n\setminus{B(x, \varepsilon)}} \frac{u(x)-u(y)}{|x-y|^{n+2s}}dy, 
    \end{array}
\end{equation}

Here $P.V.$ is a commonly used abbreviation for "in the principal value sense" (as defined by the latter equation), $x\in \mathbb{R}^n$,  $B(x, \varepsilon)$ is the ball centered in $x\in \mathbb{R}^n$ of radius $\varepsilon$ and $S(\mathbb{R}^n)$ denotes the Schwartz space.\\
 We define the Sobolev space $ H^s(\mathbb{R}^n)$ the set of functions $u$ such that they are square integrable and their fractional Laplacian $(-\Delta)^{\frac{s}{2}} u$ is also square integrable: 
  \begin{equation}
 H^s(\mathbb{R}^n) :=\big\{u : \mathbb{R}^n \rightarrow \R,  u\in L^2(\mathbb{R}^n)\  \mbox{and}\  \int_{\mathbb{R}^n \times \mathbb{R}^n}\frac{|u(x)-u(y)|^2}{|x-y|^{n+2s}}  dx dy < +\infty \big\}.
 \end{equation}
 $H^s(\mathbb{R}^n)$ is endowed with the norm defined as

  $$  ||u||_{H^s(\mathbb{R}^n)}   :=||u||_{L^2(\mathbb{R}^n)}+\bigg ( \int_{\mathbb{R}^n \times \mathbb{R}^n}\frac{|u(x)-u(y)|^2}{|x-y|^{n+2s}}  dx dy \bigg )^\frac{1}{2}.$$

We know that $H^s(\mathbb{R}^n)$ endowed with the norm $\|.\|_{H^s(\mathbb{R}^n)}$  is a Hilbert space.\\ 
Thus, $\mathbb{H}^s_0(\Omega)$ is a subspace of $H^s(\mathbb{R}^n)$ and is defined simply as in \cite{BRS} by 

$$\mathbb{H}^s_0(\Omega)= \big\{ u\in H^s(\mathbb{R}^n); u=0\mbox { a.e. sur }\Omega^c \big\}.$$

 The norm in $\mathbb{H}_0^s(\Omega)$ is defined as follows
$$\mathcal{N}(u) := ||u||_{L^2(\Omega)}+\bigg(\displaystyle\int_{\mathbb{R}^n \times \mathbb{R}^n}\frac{|u(x)-u(y)|^2}{|x-y|^{n+2s}}  dx dy \bigg )^\frac{1}{2}.$$
This norm  is equivalent to 
\begin{equation*}
||u||_{\mathbb{H}_0^s(\Omega)}= \bigg(\int_{\mathbb{R}^n \times \mathbb{R}^n}\displaystyle\frac{|u(x)-u(y)|^2}{|x-y|^{n+2s}}  dx dy \bigg )^\frac{1}{2},
\end{equation*}
see (\cite{BRS}, lemma1.28) and (\cite{FLS}, lemma 3.1).
We notice that the following identity 
$||u||_{\mathbb{H}^s_0(\Omega)}=||(-\Delta)^\frac{s}{2}u||_{L^2(\mathbb{R}^n)}$, gives the relation between the fractional Laplacian operator
$(-\Delta)^s$ and the fractional Sobolev space $\mathbb{H}^s_0(\Omega)$, see (\cite{BRS}, (4.35)).
In this present paper, $\mathbb{H}^s_0(\Omega)$ will be the functional analytic setting because the classical fractional Sobolev space approach not sufficient for studying the problem, see \cite{ABS}, \cite{BRS}, \cite{SV}.\\
We start by recalling some notations and some remarks which will be useful. First,  We will start by giving some important definitions afterwards.
Let denote by $S_s(p) := S_{s, 0}(p)$
and $S_s :=S_{s, 0}(1)$ the weightless case.
Let's denote by $u_{\varepsilon, s, a}$ an extremal function for the weightless Sobolev inequality for the fractional Laplacian operator. Let us fix $\eta>0$ such that 
\begin{equation}\label{eq eta}
B(a, {4\eta}) \subset \Omega
\end{equation}
For $x \in \mathbb{R}^n$ and $\varepsilon>0$, 
\begin{equation}\label{fonction test u}
    U_{\varepsilon, s, a}(x)=\bigg(\frac{\varepsilon}{\varepsilon^2+|x-a|^2}\bigg)^{\frac{n-2s}{2}}, 
 \end{equation}
then we set
\begin{equation}\label{fonction test u epsilon}
    u_{\varepsilon, s, a}(x)= U_{\varepsilon, s, a}(x) \Psi(x), 
\end{equation}
where $\Psi \in \C_0^\infty (\mathbb{R}^n)$ such that $0\leq\Psi\leq 1$,  $\Psi = 1$ in $B(a, \eta)$ and $\Psi=0$ in $ B(a, 2\eta)^c$,  with $\eta$ is a postive real.\\
Using \cite{BN1} and \cite{BRS}, we have
\begin{equation}\label{norme 2 lap frac u epsilon}
    \displaystyle\int_{\mathbb{R}^n}\displaystyle\int_{\mathbb{R}^n} \displaystyle\frac{|u_{\varepsilon, s, a}(x)-u_{\varepsilon, s, a}(y)|^2}{|x-y|^{n+2s}}dydx=K_{s} +  \mathcal{O}(\varepsilon ^{n-2s}),
\end{equation}
for $n\geq4$ and $s \in ]0,1[,$ or  $n=3$ and $s < \frac{3}{4}$,  the square integration of $u_{\varepsilon, s, a}$ goes as follows for any $n > 4 s$ 
\begin{equation}\label{norme 2 u epsilon}
||u_{\varepsilon, s, a}||_{L^2(\mathbb{R}^n)}^2 = K_{2, s} \varepsilon^{2s} + \mathcal{O}(\varepsilon ^{n-2s}), 
\end{equation}
the $q_s$ norm is given by
\begin{equation}\label{norme qs de u epsilon}
||u_{\varepsilon, s, a}||_{L^{q_s}(\mathbb{R}^n)}^{q_s}=K_{q_s}+ \mathcal{O}(\varepsilon^n) ,
\end{equation}
$\hbox{ with }K_{q_s}=\displaystyle\int_{\mathbb{R}^n}\frac{dy}{(1+|y|^2)^{n}}$  and $\frac{K_{s}}{K_{q_s}^\frac{2}{q_s}}=S_s.$\\

In order to present our main results, we will need to introduce the first eigenvalue $\lambda_{1,p,s}$ on $\Omega$ with zero Dirichlet boundary condition 
\begin{equation}\label{déf lambda_1,p,s}
\lambda_{1, p, s} :=\min_{\substack {u\in \mathbb{H}_0^s(\Omega)\\ u \ne 0}}   \frac{\displaystyle\int_{\mathbb{R}^n}p(x) \bigg(\displaystyle\int_{\mathbb{R}^n}\frac{|u(x)-u(y)|^2}{|x-y|^{n+2s}}dy\bigg)dx}{\displaystyle\int_{\mathbb{R}^n}|u(x)|^2dx}. 
\end{equation}

\subsection{Natural scaling(s) of the problem}
Let us define the energy 
$$E_\lambda(u) := \displaystyle\int_{\mathbb{R}^n}p(x) \bigg(\displaystyle\int_{\mathbb{R}^n}\frac{|u(x)-u(y)|^2}{|x-y|^{n+2s}}dy\bigg)dx - \lambda \displaystyle\int_{\Omega}u ^2 dx.
$$
 For the sake of clarity, using an argument in \cite{HV}, we analyze a blow up around a minimum of $p$ which we suppose  $0$,  we are led to study $E_\lambda(v)$. More precisely, let us define $v_\varepsilon$ by $v(x)=\varepsilon^{-n/q_s} v_\varepsilon(\frac{x}{\varepsilon})$ and $\Omega_\varepsilon = \varepsilon^{-1} \Omega$. Note that if  $\vert \vert  v \vert\vert_{L^{q_s}(\R^n)} =1$ then $\vert \vert  v_\varepsilon \vert\vert_{L^{q_s}(\R^n)} =1$. For $p(x)=p_0 + \kappa |x-a|^k$ in $B(a,4\eta)$ and satisfying \eqref{p-p_0 dans L1}, we have
\begin{equation*}
    \begin{array}{ll}
E_\lambda(v) & =
\displaystyle\int_{\R^n}p(x) \bigg(\displaystyle\int_{\mathbb{R}^n}\frac{|v(x)-v(y)|^2}{|x-y|^{n+2s}}dy\bigg)dx - \lambda \displaystyle\int_{\Omega}v^2 dx\\
\\
& = \displaystyle\int_{B(0,{\eta\over \varepsilon})}p(\varepsilon x) \bigg(\displaystyle\int_{\mathbb{R}^n}
\frac{|v_\varepsilon(x)-v_\varepsilon(y)|^2}{|x-y|^{n+2s}}dy\bigg)dx
\\
\\
& + \displaystyle\int_{B(0,\eta^c)}p(x) \bigg(\displaystyle\int_{\mathbb{R}^n}\frac{|v(x)-v(y)|^2}{|x-y|^{n+2s}}dy\bigg)dx - \lambda \varepsilon^{2s} \displaystyle\int_{\Omega_\varepsilon}v_\varepsilon(x)^2 dx

    \end{array}
\end{equation*}

\begin{equation*}
    \begin{array}{ll}
E_\lambda(v)
&= p_0\displaystyle\int_{\R^n}|(-\Delta)^\frac{s}{2} v(x)|^2dx +  \kappa\varepsilon^k \displaystyle
\int_{B(0,{\eta\over\varepsilon} )} |x-a|^k \displaystyle\int_{\mathbb{R}^n}
\frac{|v_\varepsilon(x)-v_\varepsilon(y)|^2}{|x-y|^{n+2s}}dydx \\
& + \displaystyle\int_{B(0,\eta)^c}p(x) \bigg(\displaystyle\int_{\mathbb{R}^n}\frac{|v(x)-v(y)|^2}{|x-y|^{n+2s}}dy\bigg)dx
- \lambda \varepsilon^{2s} \displaystyle\int_{\Omega_\varepsilon}v_\varepsilon^2 dx.
    \end{array}
\end{equation*}

Unlike to for example \cite{HV}, we have a positive term in the previous expression which is the following $\displaystyle\int_{B(0,\eta)^c}p(x) \bigg(\displaystyle\int_{\mathbb{R}^n}\frac{|v(x)-v(y)|^2}{|x-y|^{n+2s}}dy\bigg)dx$ that can explain the conditions on $k$ announced later in our result.
\\
When  $k >2s$, the term related to $\varepsilon^{2s}$ dominates the one related to $\varepsilon^k$ provided that  the quantity

\begin{equation*}\label{intégralefinie}
\int_{B(0,{\eta\over\varepsilon} )} |x-a|^k\int_{\mathbb{R}^n}\frac{|v_\varepsilon(x)-v_\varepsilon(y)|^2}{|x-y|^{n+2s}}dydx
\end{equation*}
is finite, so we find that the energy get strictly below the critical value and we expect that  $E_\lambda(v) < p_0 S_{s}$, which gives that   $S_{s, \lambda}(p)$ is achieved.\\

It is  important to  mention that, since the problem is not local, at infinity the term $|x-a|^k$ has an impact on the  above  integral  which is not finite in general for $u\in \mathbb{H}_0^s(\Omega)$,  so that we will restrict ourselves to specific values of $k$.

We expect a competition between the local character of the weight and the nonlocal operator. It turns out that the local nature of the weight wins for a class of weights. The presence of weight pushes the problem to be local. In other words, when $k$ is smaller than $n-4s$, we find that the local character dominates, on the other hand when $k$ is very large we expect that the problem does not admit solutions.

If we take $\lambda=0$  we get the following nonexistence result.
\begin{prop}
     For $n =  3$ and $s\in ]0, \frac{1}{4}]$ or $n =  4$ and $ s\in ]0,\frac{1}{2}]$ or  $n =  5$ and $s\in ]0, \frac{3}{4}]$ or $n\geq 6$ and $s\in ]0, 1[$ and for all $k \in [2,n-4s[$ we have ${S}_{s, 0}(p) = p_0 S_{s, 0}(1)$ and there is no minimizing solution for
        \begin{equation}
\begin{array}{ll}
\displaystyle\int_{\mathbb{R}^n\times \mathbb{R}^n}p(x)\frac{(u(x)-u(y))(\varphi(x)-\varphi(y))}{|x-y|^{n+2s}}dxdy    = \displaystyle\int_{\Omega} |u(x)|^{q_s-2}u(x)\varphi(x) dx,
\end{array}
\end{equation}
$\varphi \in H^s(\mathbb{R}^n).$
\end{prop}
\begin{preuve}
Let $s\in ]0, 1[$ and $2\leq k<n-4s$. Let us denote by $S_s(p) := S_{s, 0}(p)$ and $S_s :=S_{s, 0}(1)$. By Theorem \ref{inégalité p S}, please see section 2,  we write 
\begin{equation}
p_0 S_s \leq S_s(p) \leq  \displaystyle\frac{1}{||u_{\varepsilon,s,a}||^2_{L^{q_s}(\mathbb{R}^n)}}  \displaystyle\int_{\mathbb{R}^n} \displaystyle\int_{\mathbb{R}^n} p(x) \frac{|u_{\varepsilon,s,a}(x)-u_{\varepsilon,s,a}(y)|^2}{|x-y|^{n+2s}}dy dx \leq p_0 S_s+o(1),
\end{equation}
as $\varepsilon$ tends to zero.
Then we get $S_s(p)= p_0 S_s$.\\
We suppose that $S_s(p)$ is achieved by some function $u\in \mathbb{H}_0^s(\Omega)$ such that $||u||_{L^{q_s}(\mathbb{R}^n)}=1$,  then
\begin{equation}
p_0\displaystyle\int_{\mathbb{R}^n} \displaystyle\int_{\mathbb{R}^n} \frac{|u(x)-u(y)|^2}{|x-y|^{n+2s}}dy dx \leq\displaystyle\int_{\mathbb{R}^n} 
\displaystyle\int_{\mathbb{R}^n} p(x) \frac{|u(x)-u(y)|^2}{|x-y|^{n+2s}}dy dx = p_0 S_s.
\end{equation}
Therefore,  $S_s=\displaystyle\int_{\mathbb{R}^n} \displaystyle\int_{\mathbb{R}^n} \frac{|u(x)-u(y)|^2}{|x-y|^{n+2s}}dy dx,  ||u||_{L^{q_s}(\mathbb{R}^n)}=1 $ which means that $S_s$ is achieved by a function $u \in \mathbb{H}_0^s(\Omega)$ and this is absurd. See \cite{L}, see also \cite{FLS} , \cite{CT} or \cite{NPV}.
\end{preuve}

In what follows, we will concentrate on the case where  $\lambda>0$.

\subsection{Statement of the main results} 
Let us announce the main statements of this paper. We state two Theorems where in  the first, we  prove existence of minimizers of 
$S_{s, \lambda}(p)$ in the presence of a linear perturbation.

  \begin{thm}\label{théorème d'existence} 
 Let $\Omega$ be an open, bounded subset of $\mathbb{R}^n$ with continuous boundary. Let $p$ as in \eqref{def de p} and \eqref{p-p_0 dans L1}, $\lambda \in ]0, \lambda_{1, p, s}[$, $n\geq 3$ and $2\leq k<n-4s$.
The following statements hold true
\begin{enumerate}
\item
If $n\geq 6$ and $s\in ]0, 1[$, there exists $C_1 := C(n,s,k)>0$ such that for every  $\kappa \in ]0,C_1\lambda[$, $S_{s, \lambda}(p)$ is achieved.
\item 
If  $n =  3$ and $s\in ]0, \frac{1}{4}[$ or $n =  4$ and $ s\in ]0,\frac{1}{2}[$ or  $n =  5$ and $s\in ]0, \frac{3}{4}[$ 
there exists a constant $C_2 :=C( n,s,k) > 0$ such that  for every $\kappa \in ]0, C_2\lambda[$ we have  $S_{s, \lambda}(p)$ is achieved.
\end{enumerate}

  \end{thm}

The second theorem is dealing with problem \eqref{formulation faible}, we prove that it has non-ground state solutions with a subcritical pertubation  by proceeding with a min-max technique using the mountain pass Theorem see \cite{AR}.
\begin{thm}\label{prob sans minimisation}
Let $\Omega$ be an open, bounded subset of $\mathbb{R}^n$ with continuous boundary. Let $p$ be as defined in \eqref{def de p} and \eqref{p-p_0 dans L1}, we suppose also that $p$ is bounded. Let  Suppose that   $k \in [2,n-4s[$.  
If $2<q<q_s$ then for all $\lambda>0$, we have
\begin{enumerate}
\item

If $n\geq 6$ then for every $s \in ]0,1[$  problem \eqref{formulation faible} has a nontrivial solution 
$u \in \mathbb{H}_0^s(\Omega)$.
\item
If  $n =  3$ and $s\in ]0, \frac{1}{4}[$ or $n =  4$ and $ s\in ]0,\frac{1}{2}[$ or  $n =  5$ and $s\in ]0, \frac{3}{4}[$ then problem \eqref{formulation faible} has a nontrivial solution 
$u \in \mathbb{H}_0^s(\Omega)$.
\end{enumerate}
In the case where  $q= 2$ then 
\begin{enumerate}
\item
 If $n\geq 6$ and $s\in ]0, 1[$, there exists $C_1 := C(n,s,k)>0$ such that for every $\lambda \in ]0, \lambda_{1,p,s}[$ and for every $\kappa \in ]0,C_1\lambda[$, problem \eqref{formulation faible} has a nontrivial solution $u \in \mathbb{H}_0^s(\Omega)$.
\item
While if  $n =  3$ and $s\in ]0, \frac{1}{4}[$ or $n =  4$ and $ s\in ]0,\frac{1}{2}[$ or  $n =  5$ and $s\in ]0, \frac{3}{4}[$
there exists a constant $C_2 :=C( n,s,k) > 0$ such that for every $\lambda \in ]0, \lambda_{1,p,s}[$  and that  for every $\kappa \in ]0, C_2\lambda[$ we have  problem \eqref{formulation faible} has a nontrivial solution $u \in \mathbb{H}_0^s(\Omega)$.
\end{enumerate}
\end{thm}

\subsection{Structure of the paper}
The paper is structured as follows. The next section $\S2$,  proves the a-priori estimate $S_{s, \lambda}(p)<p_0 S_s$.
Theorem \ref{théorème d'existence} which is the first main result of this paper, is proved by mathematical adequate technique using the previous section in order to get existence of minimizing solutions to $S_{s,\lambda}(p)$. 
In section $\S3$ we carry out non-ground state solutions using the minimax technique and proving the mountain pass Theorem in the general case.
We investigate also  a  subcritical approximation and we adopt the strategy used  from \cite{BN1} to prove Theorem \ref{prob sans minimisation} by choosing a suitable test function.
\section{Existence of minimizers}
First of all,  let us prove that the infimum $S_{s, \lambda}(p)$ does exist for every $\lambda \in \mathbb{R}$. In fact, thanks to Sobolev inequalities we have
\begin{equation}
    \begin{array}{ll}
    \displaystyle\int_{\mathbb{R}^n}p(x) \bigg(\displaystyle\int_{\mathbb{R}^n}\frac{|u(x)-u(y)|^2}{|x-y|^{n+2s}}dy\bigg)dx-\lambda \displaystyle\int_\Omega |u(x)|^2dx 
    & \geq C_1 ||u||_{L^{q_s}(\mathbb{R}^n)} -\lambda C_2 ||u||_{L^{q_s}(\mathbb{R}^n)},
    \end{array}
\end{equation}
where $C_1, C_2$ are positive constants.
Then we get
\begin{equation}
  S_{s,\lambda}(p)  \geq   C_1-\lambda C_2
\end{equation}
Thus, the infimum $S_{s,\lambda}(p)$ exists.\\
Here, we used the fact that $ N(u):=\displaystyle\int_{\mathbb{R}^n}p(x) \bigg(\displaystyle\int_{\mathbb{R}^n}\frac{|u(x)-u(y)|^2}{|x-y|^{n+2s}}dy\bigg)dx$ is a norm on $\mathbb{H}_0^s(\Omega)$.

The following Theorem plays a crucial role to prove existence of solutions,  it is an adaptation of an original argument due to \cite{BN1} in the context of Yamabe's conjecture.
\begin{prop}\label{prop fondamentale}
Let $s\in]0, 1[$ and let a weight $p$ satisfying \eqref{def de p}.
 If $S_{s, \lambda}(p)< p_0 S_s$,  then $S_{s, \lambda}(p)$ is achieved.
\end{prop}
\begin{preuve}
Let us recall the expression of $S_{s, \lambda}(p)$:
 \begin{equation*}
 S_{s, \lambda}(p)=\inf_{\substack{u\in \mathbb{H}_0^s(\Omega)\\ ||u||_{L^{q_s}(\R^n)}=1}}\displaystyle\int_{\R^n}p(x) \bigg(\displaystyle\int_{\R^n}\frac{|u(x)-u(y)|^2}{|x-y|^{n+2s}}dy\bigg)dx-\lambda \displaystyle\int_\Omega |u(x)|^2dx.
 \end{equation*}
Let $(u_j)$ a minimizing sequence of $S_{s, \lambda}(p)$,  then $||u_j||_{L^{q_s}(\R^n)}=1$
and
\begin{equation*}
\displaystyle\int_{\R^n}p(x) \displaystyle\int_{\R^n}\frac{|u_j(x)-u_j(y)|^2}{|x-y|^{n+2s}}dydx-\lambda ||u_j||_2^2=S_{s, \lambda}(p)+o(1),
\end{equation*}
as j tends to $+\infty$.\\
We define a scalar product on $\mathbb{H}_0^s(\Omega)$ by
\begin{equation}\label{produit scalaire p}
<u, v>_p=\displaystyle\int_{\R^n} p(x) \displaystyle\int_{\R^n}\frac{(u(x)-u(y)) (v(x)-v(y))}{|x-y|^{n+2s}} dy dx,
\end{equation}
for all $u, v \in \mathbb{H}_0^s(\Omega)$.
It's obvious that the norm associated with this scalar product is equivalent to the ordinary norm over $\mathbb{H}_0^s(\Omega).$

Since $(u_j)$ is bounded in $\mathbb{H}_0^s(\Omega)$,  we extract a subsequence still denoted by $(u_j)$ such that, $(u_j)$ tends weakly to $u$ in $\mathbb{H}_0^s(\Omega)$ (since $\mathbb{H}_0^s(\Omega)$ is reflexive space). Then
$(u_j)$ tends strongly to $u$ in $L^2(\Omega)$,  and
$(u_j)$ tends to $u$ almost everywhere in $\Omega$, with
\begin{equation*}\label{majoration u par 1}
||u||_{L^{q_s}(\R^n)} \leq 1.
\end{equation*}
Using the defintion of $S_s$ we have
$$\displaystyle\int_{\R^n}p(x)
\displaystyle\int_{\R^n}\frac{|u_j(x)-u_j(y)|^2}{|x-y|^{n+2s}}dydx \geq p_0 S_s, $$

then
$$\displaystyle\int_{\R^n}p(x)
\displaystyle\int_{\R^n}\frac{|u_j(x)-u_j(y)|^2}{|x-y|^{n+2s}}dydx-S_{s, \lambda}(p)+o(1) \geq p_0 S_s-S_{s, \lambda}(p)+o(1), $$
as $j$ tends to $+\infty,$
therefore
\begin{equation*}
\lambda ||u||_{L^{2}(\R^n)}^2 \geq p_0 S_s - S_{s, \lambda}(p)>0.
\end{equation*}
Therefore $u \neq 0$.

Let's take $v_j:=u_j -u$ then $(v_j)$ tends weakly to 0 in $\mathbb{H}_0^s(\Omega)$ and $(v_j)$ tends to $0$ almost everywhere in $\Omega$. Using the definition of $S_{s,\lambda}(p)$,  we obtain
\begin{equation}\label{b}
\begin{array}{ll}
\displaystyle\int_{\R^n}p(x)\bigg(
\displaystyle\int_{\R^n}\frac{|u(x)-u(y)|^2}{|x-y|^{n+2s}}dy\bigg)dx+\displaystyle\int_{\R^n}p(x)
\displaystyle\int_{\R^n}\frac{|v_j(x)-v_j(y)|^2}{|x-y|^{n+2s}}dydx \\
\\
-\lambda ||u||_{L^{2}(\R^n)}^2  +  o(1)=  S_{s, \lambda}(p)+o(1),
\end{array}
\end{equation}
as $j$ tends to $+\infty$.\\

On the other hand,  Since $(v_j)$ is bounded in $L^{q_s}(\Omega)$ and $(v_j)$ tends to 0 almost everywhere in $\Omega$, we deduce from a result of Brezis-Lieb  that
\begin{equation*}
1= ||u+v_j||_{L^{q_s}(\R^n)}^{q_s} = \|u||_{L^{q_s}(\R^n)}^{q_s}+\|v_j||_{L^{q_s}(\R^n)}^{q_s}+o(1),
\end{equation*}
as $j$ tends to $+\infty.$
Therefore
\begin{equation*}
1  \leq \|u||_{L^{q_s}(\R^n)}^{2}+\|v_j||_{L^{q_s}(\R^n)}^{2}+o(1),
\end{equation*}
as $j$ tends to $+\infty.$
Denoting by $\Tilde{v}_j :=\frac{v_j}{||v_j||_{L^{q_s}(\R^n)}}$. Since we have
\begin{equation*}
p_0 S_s\leq \displaystyle\displaystyle\int_{\R^n}p(x)
\displaystyle\int_{\R^n}\displaystyle\frac{|\Tilde{v}_j(x)-\Tilde{v}_j(y)|^2}{|x-y|^{n+2s}}dydx,
\end{equation*}

we obtain
\begin{equation*}\label{a}
1\leq\|u||_{L^{q_s}(\R^n)}^2+\frac{1}{p_0 S_s}\displaystyle\displaystyle\int_{\R^n}p(x)
\displaystyle\int_{\R^n}\frac{|v_j(x)-v_j(y)|^2}{|x-y|^{n+2s}}dydx+o(1),
\end{equation*}
as $j$ tends to $+\infty.$
In the case where  $S_{s, \lambda}(p)\geq0$,  we deduce from \eqref{a} that
\begin{equation}\label{c}
 S_{s, \lambda}(p)\leq S_{s, \lambda}(p)\|u||_{L^{q_s}(\R^n)}^2+\frac{S_{s, \lambda}(p)}{p_0 S_s}\displaystyle\displaystyle\int_{\R^n}p(x)
\displaystyle\int_{\R^n}\frac{|v_j(x)-v_j(y)|^2}{|x-y|^{n+2s}}dydx+o(1),
\end{equation}
as $j$ tends to $+\infty.$

Combining \eqref{b} and \eqref{c},  we obtain
\begin{equation*}
\begin{array}{ll}
\displaystyle\displaystyle\int_{\R^n}p(x)\bigg(
\displaystyle\int_{\R^n}\frac{|u(x)-u(y)|^2}{|x-y|^{n+2s}}dy\bigg)dx-\lambda \displaystyle\int_\Omega|u(x)|^2dx+\\\displaystyle\displaystyle\int_{\R^n}p(x)
\displaystyle\int_{\R^n}\frac{|v_j(x)-v_j(y)|^2}{|x-y|^{n+2s}}dydx&\\
\\ \leq S_{s, \lambda}(p)||u||_{L^{q_s}(\R^n)}^2+\displaystyle\frac{S_{s, \lambda}(p)}{p_0 S_s}\displaystyle\displaystyle\int_{\R^n}p(x)
\displaystyle\int_{\R^n}\frac{|v_j(x)-v_j(y)|^2}{|x-y|^{n+2s}}dydx+o(1),&
\end{array}
\end{equation*}
as $j$ tends to $+\infty.$

Thus,
\begin{equation*}
    \begin{array}{ll}
\displaystyle\displaystyle\int_{\R^n}p(x)\bigg(
\displaystyle\int_{\R^n}\frac{|u(x)-u(y)|^2}{|x-y|^{n+2s}}dy\bigg)dx-\lambda \displaystyle\int_\Omega|u(x)|^2dx  +o(1) \leq &\\
S_{s, \lambda}(p)||u||_{L^{q_s}(\R^n)}^2+\bigg[\displaystyle\frac{S_{s, \lambda}(p)}{p_0 S_s}-1\bigg]\displaystyle\displaystyle\int_{\R^n}p(x)
\displaystyle\int_{\R^n}\frac{|v_j(x)-v_j(y)|^2}{|x-y|^{n+2s}}dydx,
\end{array}
\end{equation*}
as $j$ tends to $+\infty.$
Since $S_{s, \lambda}(p)< p_0 S_s$,  we deduce
\begin{equation}\label{d}
\displaystyle\displaystyle\int_{\R^n}p(x)\bigg(
\displaystyle\int_{\R^n}\frac{|u(x)-u(y)|^2}{|x-y|^{n+2s}}dy\bigg)dx-\lambda \displaystyle\int_\Omega|u(x)|^2dx\leq S_{s, \lambda}(p)||u||_{L^{q_s}(\R^n)}^2  .
\end{equation}
This means that $u$ is a minimum of $S_{s, \lambda}(p)$.
 \end{preuve}
\begin{remarque}
In the case where the weight $p$ and $\lambda$ are such that $-\infty<S_{s, \lambda}(p)\leq 0$,  we prove that the infimum $S_{s, \lambda}(p)$ is achieved. In fact,  as in \eqref{majoration u par 1} we have $||u||_{L^{q_s}(\mathbb{R}^n)}\leq 1$,  then
$$S_{s, \lambda}(p)\leq S_{s, \lambda}(p)||u||_{L^{q_s}(\mathbb{R}^n)}^2.$$
Again,  we deduce \eqref{d} from \eqref{b}.
\end{remarque}
 
\begin{remarque} \label{solution positive}
Thanks to the following inequality $\big| |u_j(x)| - |u_j(y)|  \big| \leq |u_j(x) - u_j(y)|$,  if $(u_j)$ is a minimising sequence of $S_{s, \lambda}(p)$,  $(|u_j|)$ is also a minimising sequence. Therefore,  we can take a positive minimizing sequence of $S_{s, \lambda}(p)$.

Clearly if $S_{s,\lambda}(p)$ is achieved by some function $u$, then it is also achieved by $|u|$. Therefore, there exists a positive solution for the infimum $S_{s,\lambda}(p).$
\end{remarque}
\subsection{ A priori estimate on $S_{s,\lambda}(p)$ }
We need to prove the following theorem and the techniques that it uses since it constitutes an official key in order to apply Proposition \ref{prop fondamentale}.
\begin{thm} \label{inégalité p S}
Let $n,s$ such that $n =  3$ and $s\in ]0, \frac{1}{4}[$ or $n =  4$ and $ s\in ]0,\frac{1}{2}[$ or  $n =  5$ and $s\in ]0, \frac{3}{4}[$ or $n\geq 6$ and $s\in ]0, 1[$. Let $k \in [2,n-4s[$. Then the following estimate holds true
 \begin{equation}\label{inégalité clé}
 \displaystyle\int_{\mathbb{R}^n}\displaystyle\int_ {\mathbb{R}^n} p(x)  \frac{|u_{\varepsilon, s, a}(x)-u_{\varepsilon, s, a}(y)|^2}{|x-y|^{n+2s}}dx dy\leq p_0 K_s +\kappa C\varepsilon^{2s}+ \mathcal{O}(\varepsilon^{n-2s}) + \mathcal{O}(\varepsilon^{k+2s}), 
 \end{equation}
 as $\varepsilon$ tends to zero, where  $K_s$ is defined in \eqref{norme 2 lap frac u epsilon} and $C$ is a positive constant depending on $k$, $s$ and the dimension $n$.
 \end{thm}
 \begin{preuve}
 First of all and without loss of generality, we assume that $a=0$ and we note $U_{\varepsilon,s,0}= U_{\varepsilon,s}$  $u_{\varepsilon,s,0}= u_{\varepsilon,s}$.\\
 In the sequel, we will assume in that $n$ and $s$ are such that
 \begin{equation}\label{conditions n s}
     n =  3 \hbox{ and } s\in ]0, \frac{1}{4}[\\
     \\ \hbox{ or } n =  4 \hbox{ and } s\in ]0,\frac{1}{2}[ \\
     \\ \hbox{ or }  n =  5 \hbox{ and } s\in ]0, \frac{3}{4}[ \\
     \\ \hbox{ or } n\geq 6 \hbox{ and } s\in ]0, 1[.
 \end{equation}
 \\
 The proof makes use of the following estimates and is a bit complicated than the one for similar results in the case of the fractional Laplacian without weight as it is in \cite{BRS} that we will inspire from it to accomplish the proof.\\
Using \eqref{norme 2 lap frac u epsilon}, we have
 \begin{equation}\label{10}
     \begin{array}{ll}
\displaystyle\int_{{\mathbb{R}^n}\times{\mathbb{R}^n}} p(x)  \frac{|u_{\varepsilon, s, a}(x)-u_{\varepsilon, s, a}(y)|^2}{|x-y|^{n+2s}}dx dy
\leq p_0 K_s + \displaystyle\int_{{\mathbb{R}^n}\times{\mathbb{R}^n}} (p(x)-p_0)  \frac{|u_{\varepsilon, s, a}(x)-u_{\varepsilon, s, a}(y)|^2}{|x-y|^{n+2s}}dx dy
     \end{array}
 \end{equation}
 Our aim is to prove the following inequality
\begin{equation}
    \displaystyle\int_{\mathbb{R}^n}\displaystyle\int_ {\mathbb{R}^n} (p(x)-p_0)  \frac{|u_{\varepsilon, s, a}(x)-u_{\varepsilon, s, a}(y)|^2}{|x-y|^{n+2s}}dx dy \leq \kappa C\varepsilon^{2s}+ \mathcal{O}(\varepsilon^{n-2s}) + \mathcal{O}(\varepsilon^{k+2s})
\end{equation}
 We easily see that for $\rho >0$ and $x\in B(0, \rho)^c$,  then
 \begin{equation}\label{majoration fonction test u}
 u_{\varepsilon,s}(x)\leq U_{\varepsilon, s}(x)\leq C \varepsilon^{\frac{n-2s}{2}}, 
 \end{equation}
 for any $\varepsilon>0$ and for some positive constant $C$,  possibly depending on $\eta$, $\rho$, $s$ and the dimension $n$.
\\
We mention that the following assertions hold true
 \\
 (a) For any $x\in\mathbb{R}^n$ and $y \in B(0, \eta)^c$,  with $|x-y|\leq\frac{\eta}{2}, $
 \begin{equation}\label{a)}
     |u_{\varepsilon, s}(x)-u_{\varepsilon,s}(y)| \leq C \varepsilon^\frac{n-2s}{2}|x-y|.
 \end{equation}
 \\
 (b) For any $x, y \in B(0, \eta)^c, $
 \begin{equation}\label{b)}
     |u_{\varepsilon,s}(x)-u_{\varepsilon,s}(y)| \leq C \varepsilon^\frac{n-2s}{2}\min\{1, |x-y|\}, 
 \end{equation}
for any $\varepsilon>0$ and for some positive constant $C$,  possibly depending on $\eta$, $\rho$, $s$ and $n$
\\
 Let us begin the proof of the Theorem. We introduce the notations
 \begin{equation}
 \mathbb{D} := \{(x, y)\in \mathbb{R}^n \times \mathbb{R}^n : x\in B(0, \eta) ,  y\in B(0, \eta)^c ,  |x-y|>\frac{\eta}{2}   \}
 \end{equation}
 and
 \begin{equation}
 \mathbb{E} :=\{(x, y)\in \mathbb{R}^n \times \mathbb{R}^n : x\in B(0, \eta) ,  y\in B(0, \eta)^c ,  |x-y|\leq\frac{\eta}{2}   \},  
 \end{equation}
 where $\eta$ is as in \eqref{eq eta}.\\
By \eqref{fonction test u epsilon}, and using the properties of $p$ described in \eqref{def de p} and \eqref{p-p_0 dans L1} we have that
 \begin{equation}\label{4 intégrales}
   \begin{array}{ll}
   \displaystyle\int_{\mathbb{R}^n}\displaystyle\int_ {\mathbb{R}^n} (p(x)-p_0)\displaystyle\frac{|u_{\varepsilon,s}(x)-u_{\varepsilon,s}(y)|^2}{|x-y|^{n+2s}}dy dx 
   \\
   \\
   \leq \kappa \displaystyle\int_{B(0, \eta)}\displaystyle\int_ {B(0, \eta)}|x|^k\frac{|U_{\varepsilon, s}(x)-U_{\varepsilon, s}(y)|^2}{|x-y|^{n+2s}}dy dx \\
   \\
    + 2 \kappa \displaystyle\int_{\mathbb{D}} |x|^k\frac{|u_{\varepsilon,s}(x)-u_{\varepsilon,s}(y)|^2}{|x-y|^{n+2s}}dy dx \\
   \\
   +2 \kappa \displaystyle\int_ {\mathbb{E}} |x|^k\frac{|u_{\varepsilon,s}(x)-u_{\varepsilon,s}(y)|^2}{|x-y|^{n+2s}}dy dx \\
   \\
     +  \displaystyle\int_{B(0, \eta)^c} \displaystyle\int_ {B(0, \eta)^c}  (p(x)-p_0) \frac{|u_{\varepsilon,s}(x)-u_{\varepsilon,s}(y)|^2}{|x-y|^{n+2s}}dy dx
   \end{array}
 \end{equation}
 
We start by treating the first term in the right hand side.
For $n$ and $s$ as described in \eqref{conditions n s},  $2\leq k<n-4s$ and $\varepsilon>0$,  let us denote by 
\begin{equation}\label{A_s,k,epsilon}
A_{s, k, \varepsilon} :=\displaystyle\int_{|x|\leq \eta} |x|^k \bigg ( \displaystyle\int_{|y|\leq \eta}\frac{|U_{\varepsilon, s}(x)-U_{\varepsilon, s}(y)|^2}{|x-y|^{n+2s}}dy\Bigg ) dx.  \\
\end{equation}

\begin{prop}\label{As,k}
For $n$ and $s$ as described in \eqref{conditions n s}, $2\leq k<n-4s$ and $\varepsilon>0$,  we have
\begin{equation}\label{estimation A}
A_{s, k, \varepsilon} \leq C \varepsilon^{2s},
\end{equation}
   where $C$ is a positive constant depending on $s,k$ and the dimension $n$.
\end{prop}
\begin{preuve} 
Performing a change of variables in \eqref{A_s,k,epsilon},  we obtain
\begin{equation}\label{A epsilon}
    \begin{array}{ll}
      A_{s, k, \varepsilon}   & =\varepsilon^k\displaystyle\int_{|x|\leq\frac{\eta}{\varepsilon}}  |x|^k \displaystyle\int_{|y|\leq\frac{\eta}{\varepsilon}}\frac{1}{|x-y|^{n+2s}}\Bigg |\frac{1}{(1+|x|^2)^\frac{n-2s}{2}}-\frac{1}{(1+|y|^2)^\frac{n-2s}{2}} \Bigg |^2 dy dx  \\
         & = \varepsilon^k\displaystyle\int_{|x|\leq\frac{\eta}{\varepsilon}}\displaystyle\int_{|y|\leq\frac{\eta}{\varepsilon}} \bigg |\frac{|x|^\frac{k}{2}}{(1+|x|^2)^\frac{n-2s}{2}} - \frac{|x|^\frac{k}{2}}{(1+|y|^2)^\frac{n-2s}{2}}\bigg |^2\frac{dy dx}{|x-y|^{n+2s}}, 
    \end{array}
\end{equation}
on the other hand,  we have
\begin{equation}\label{inégalité A}
\Bigg |\frac{|x|^\frac{k}{2}}{(1+|x|^2)^\frac{n-2s}{2}}-\frac{|x|^\frac{k}{2}}{(1+|y|^2)^\frac{n-2s}{2}} \Bigg |^2\leq 2 \Bigg( |f_k(x)-f_k(y)|^2+ \bigg |\frac{|y|^\frac{k}{2}-|x|^\frac{k}{2}}{(1+|y|^2)^\frac{n-2s}{2}} \bigg |^2 \Bigg) 
\end{equation}
with $f_k(t) :=\displaystyle\frac{|t|^\frac{k}{2}}{(1+|t|^2)^\frac{n-2s}{2}},  t\in \mathbb{R}^n$.\\
Thus,
\begin{equation}\label{inégalité clés}
 A_{s, k, \varepsilon} \leq  2\varepsilon^k\displaystyle\int_{|x|\leq\frac{\eta}{\varepsilon}}\displaystyle\int_{|y|\leq\frac{\eta}{\varepsilon}} \Bigg( |f_k(x)-f_k(y)|^2+ \bigg |\displaystyle\frac{|y|^\frac{k}{2}-|x|^\frac{k}{2}}{(1+|y|^2)^\frac{n-2s}{2}} \bigg |^2 \Bigg) \frac{dy dx}{|x-y|^{n+2s}}
 \end{equation}
As a consequence,  in order to obtain the inequality \eqref{estimation A},  we will prove these two following assertions
\begin{equation}\label{1}
    \displaystyle\int_{|x|\leq\frac{\eta}{\varepsilon}} \displaystyle\int_{|y|\leq \frac{\eta}{\varepsilon}}\displaystyle\frac{|f_k(x)-f_k(y)|^2}{|x-y|^{n+2s}}dy dx \leq \bar{\tilde{C}}, 
\end{equation}
where $\bar{\tilde{C}}$ is a positive constant possibly depending on $k, s, \eta$ and $n$ 
and
\begin{equation}\label{2}
\displaystyle\int_{|x|\leq\frac{\eta}{\varepsilon}} \displaystyle\int_{|y|\leq\frac{\eta}{\varepsilon}}\displaystyle\frac{\bigg| |x|^\frac{k}{2}-|y|^\frac{k}{2} \bigg|^2}{(1+|y|^2)^{n-2s} |x-y|^{n+2s}} dy dx \leq C' + C'' \varepsilon^{-k+2s},
\end{equation}
where $C'$ and $C''$ are two positive constants possibly depending on $s, k, \eta$ and $n$. 

We start proving \eqref{1}. For $A$ a subset domain of $\mathbb{R}^n$,  we define the function $1_A$ by $\forall x \in \mathbb{R}^n$,  $1_A(x)= 0$ if $x\notin A$ and $1_A(x)=1$ if $x\in A$.\\
Let $\gamma > 0$. at first we take the case where $|x-y|\geq \gamma$.
\begin {equation}
     \begin{array}{ll}
      \displaystyle\int_{\mathbb{R}^n} \displaystyle\int_{|x-y| \geq \gamma}\frac{|f_k(x)-f_k(y)|^2}{|x-y|^{n+2s}}dy dx
      &\\
      \\ \leq 2 \Bigg(\displaystyle\int_{\mathbb{R}^n} \displaystyle\int_{|x-y| \geq \gamma}\frac{|x|^k}{(1+|x|^2)^{n-2s} |x-y|^{n+2s}} dy dx+\displaystyle\int_{\mathbb{R}^n} \displaystyle\int_{|x-y| \geq \gamma}\frac{|y|^k}{(1+|y|^2)^{n-2s} |x-y|^{n+2s}}  \Bigg)& \\
       \\
       =2 \Bigg(\displaystyle\int_{{\mathbb{R}^n} \times{\mathbb{R}^n}}\frac{|x|^k}{(1+|x|^2)^{n-2s}}\frac{1_{B(0, \gamma)^c}(x-y)}{|x-y|^{n+2s}}dy dx + \displaystyle\int_{{\mathbb{R}^n} \times{\mathbb{R}^n}}\frac{|y|^k}{(1+|y|^2)^{n-2s}}\frac{1_{B(0, \gamma)^c}(x-y)}{|x-y|^{n+2s}}dy dx \Bigg)& \\
       \\
        =4\displaystyle\int_{\mathbb{R}^n} (g_k*h)(y)dy, \\
        \end{array}
        \end{equation}
\\
where $g_k$ and $h$ are functions defined by,  for all $t \in \mathbb{R}^n$, 
$$g_k(t) :=f_k^2(t)=\frac{|t|^k}{(1+|t|^2)^{n-2s}},\   \hbox{and} \ h(t) :=\frac{1_{B(0, \gamma)^c}(t)}{|t|^{n+2s}}.$$
Provided that $k<n-4s$ (for $n$ and $s$ as in \eqref{conditions n s}),  we have $g_k\in L^1(\mathbb{R}^n)$. Since $s>0$,  $h\in L^1(\mathbb{R}^n)$.
This implies $g_k*h \in L^1(\mathbb{R}^n)$ and we have
$$||g_k*h||_{L^1(\mathbb{R}^n)}\leq ||g_k||_{L^1(\mathbb{R}^n)} ||h||_{L^1(\mathbb{R}^n)} < +\infty, $$
thus, 
\begin{equation}\label{sup à gamma}
\displaystyle\int_{\mathbb{R}^n} \displaystyle\int_{|x-y| \geq \gamma}\frac{|f_k(x)-f_k(y)|^2}{|x-y|^{n+2s}}dy dx < +\infty,
\end{equation}
\\
Now, let's take the case where $|x-y| < \gamma$.
We have 
\begin{equation}\label{taylor integrale}
f_k(x)-f_k(y)= \displaystyle\int_0^1 Df_k(y+t(x-y)) (x-y) dt
\end{equation}
with  $f_k(z) :=\displaystyle\frac{|z|^\frac{k}{2}}{(1+|z|)^\frac{n-2s}{2}},  z \in \mathbb{R}^n$,
and
$$Df_k(z)=\displaystyle\frac{|z|^{\frac{k}{2}-2}z}{(1+ |z|^2)^\frac{n-2s}{2}}\bigg(\frac{k}{2}-(n-2s)\displaystyle\frac{|z|^2}{(1+|z|^2)}\bigg)$$

Let $R>0$ such that $R\geq 2\gamma$ and let $|y| \geq R$.\\
For the sake of clarity, we note that since we have $|x-y|<\gamma$ and $|y|\geq R$, so we can suppose also that $|x| \geq R$.\\
Since $t \in [0,1]$  and $|x-y|< \gamma$, then
\begin{equation}
    \begin{array}{ll}
 |y+t(x-y)| & \geq |y|-t|x-y| \\
 & \geq |y| - \gamma \\
 & \geq |y| (1-\frac{\gamma}{|y|})\\
 & \geq \frac{1}{2} |y|
    \end{array}
\end{equation}
So we have
\begin{equation}\label{majoration différentielle}
\begin{array}{ll}
|Df_k(y+t(x-y))(x-y)| & \leq C \displaystyle\frac{|y+t(x-y)|^{\frac{k}{2}-1}|x-y|}{(1+|y+t(x-y)|^2)^\frac{n-2s}{2}} \\
 \\
& \leq C \displaystyle\frac{|x-y|}{|y+t(x-y)|^{n-2s-\frac{k}{2}+1}}\\
\\
& \leq C \displaystyle\frac{|x-y|}{|y|^{n-2s-\frac{k}{2}+1}},
\end{array}
\end{equation}
where $C$ is a constant.\\
Thus by \eqref{majoration différentielle} and \eqref{taylor integrale}
\begin{equation}
|f_k(x)-f_k(y)| \leq |Df_k(y+t(x-y))(x-y)| \leq C \displaystyle\frac{|x-y|}{|y|^{n-2s-\frac{k}{2}+1}}
\end{equation}
Since we have $k<n-4s$ (for $n$ and $s$ as in \eqref{conditions n s}) and $s<1$, we have
\begin{equation}\label{première estimation}
    \begin{array}{ll}
\displaystyle\int_{|y| \geq R} \displaystyle\int_{ |x|\geq R, |x-y| < \gamma}\displaystyle\frac{|f_k(x)-f_k(y)|^2}{|x-y|^{n+2s}}dx dy \\
\\
\leq C \displaystyle\int_{|y| \geq R} \displaystyle\int_{ |x|\geq R, |x-y| < \gamma} \displaystyle\frac{1}{|x-y|^{n+2s-2}|y|^{2n-4s-k+2}} dx dy \\
\\
 \leq C \displaystyle\int_{|y| \geq R} \displaystyle\frac{1}{|y|^{2n-4s-k+2}} \bigg(\displaystyle\int_{ |x|\geq R, |x-y| < \gamma} \displaystyle\frac{1}{|x-y|^{n+2s-2}}  dx \bigg) dy\\
\\
 \leq C \displaystyle\int_{|y| \geq R} \displaystyle\frac{1}{|y|^{2n-4s-k+2}} \bigg(\displaystyle\int_{ |y+z|\geq R, |z| < \gamma} \displaystyle\frac{1}{|z|^{n+2s-2}}  dz \bigg) dy\\
\\
 \leq C \displaystyle\int_{|y| \geq R} \displaystyle\frac{1}{|y|^{2n-4s-k+2}} \bigg(\displaystyle\int_{ |z| < \gamma} \displaystyle\frac{1}{|z|^{n+2s-2}}  dz \bigg) dy\\
\\
  \leq C \displaystyle\int_{|y| \geq R} \displaystyle\frac{dy}{|y|^{2n-4s-k+2}}   \displaystyle\int_{ |z| < \gamma} \displaystyle\frac{dz}{|z|^{n+2s-2}}   \\
\\
 < +\infty
    \end{array}
\end{equation}
Now let $|y|<R$ and so we take $|x|<R$ since $|x-y|<\gamma$. \\
Since for all $t\in[0,1]$, $|y+t(x-y|<2R$ and $\displaystyle\frac{1}{(1+|y+t(x-y)|^2)^\frac{n-2s}{2}} \leq 1,$ we have
\begin{equation}
\begin{array}{ll}
 |Df_k(y+t(x-y))(x-y)| & \leq C \displaystyle\frac{|y+t(x-y)|^{\frac{k}{2}-1}|x-y|}{(1+|y+t(x-y)|^2)^\frac{n-2s}{2}} \\
 \\
 & \leq C' |x-y|,
\end{array}
\end{equation}
where $C'$ is a positive constant.\\
Since $s<1$, we have
\begin{equation}\label{deuxième estimation}
    \begin{array}{ll}
\displaystyle\int_{|y|<R,|x|<R, |x-y| < \gamma}\frac{|f_k(x)-f_k(y)|^2}{|x-y|^{n+2s}}dy dx & \leq C' \displaystyle\int_{|y|<R,|x|<R, |x-y| < \gamma} \displaystyle\frac{|x-y|^2}{|x-y|^{n+2s}}dx dy \\
\\
& \leq C' \displaystyle\int_{|y|<R} \bigg(\displaystyle\int_{|x|<R, |x-y| < \gamma} \displaystyle\frac{dx}{|x-y|^{n+2s-2}}\bigg) dy \\
\\
& \leq C' \displaystyle\int_{|y|< R} \bigg(\displaystyle\int_{|y+z|<R, |z|< \gamma}\displaystyle\frac{dz}{|z|^{n+2s-2}}\bigg) dy  \\
\\
& \leq C' \displaystyle\int_{|y|< R} \bigg(\displaystyle\int_{|z|< \gamma}\displaystyle\frac{dz}{|z|^{n+2s-2}}\bigg) dy \\
\\
& < +\infty.
    \end{array}
\end{equation}

Finally by \eqref{première estimation} and \eqref{deuxième estimation} we get
\begin{equation}\label{inf à gamma}
\begin{array}{ll}
\displaystyle\int_{\mathbb{R}^n} \displaystyle\int_{|x-y| < \gamma}\frac{|f_k(x)-f_k(y)|^2}{|x-y|^{n+2s}}dy dx & \leq  \displaystyle\int_{|y|<R,|x|<R, |x-y| < \gamma}\frac{|f_k(x)-f_k(y)|^2}{|x-y|^{n+2s}}dy dx \\
\\
& + \displaystyle\int_{|y| \geq R, |x|\geq R, |x-y| < \gamma}\frac{|f_k(x)-f_k(y)|^2}{|x-y|^{n+2s}}dy dx \\
\\
& < +\infty
\end{array}
\end{equation}

Therefore, using \eqref{sup à gamma} and \eqref{inf à gamma},
\begin{equation}\label{f_k finie}
    \begin{array}{ll}
        \displaystyle\int_{|x|\leq\frac{\eta}{\varepsilon}} \displaystyle\int_{|y|\leq\frac{\eta}{\varepsilon}}\frac{|f_k(x)-f_k(y)|^2}{|x-y|^{n+2s}}dy dx & \leq \displaystyle\int_{\mathbb{R}^n} \displaystyle\int_{\mathbb{R}^n}\frac{|f_k(x)-f_k(y)|^2}{|x-y|^{n+2s}}dy dx \\
        \\
        & = \displaystyle\int_{\mathbb{R}^n} \displaystyle\int_{|x-y|<\gamma}\frac{|f_k(x)-f_k(y)|^2}{|x-y|^{n+2s}}dy dx \\
        \\
        & + \displaystyle\int_{\mathbb{R}^n} \displaystyle\int_{|x-y|\geq \gamma}\frac{|f_k(x)-f_k(y)|^2}{|x-y|^{n+2s}}dy dx \\
       \\
        & \leq \bar{\tilde{C}} < +\infty,
    \end{array}
\end{equation}
where $\bar{\tilde{C}}$ is a positive constant depending on $s,k$ and the dimension $n$. Thus, \eqref{1} holds true.
\\
Now, in order to prove \eqref{2},  we will need to prove the following lemma.

\begin{lemme}\label{lemme de première intégrale} Let $\gamma,  R \in \R^*_+$ and $s \in ]0, 1[$.
Let $k\geq 2$ and let the set $\mathcal{A}$ defined by $\mathcal{A}:= \{ x, y \in \mathbb{R}^n,  |x-y|<\gamma,  |x|\leq R,  |y|\leq R \}$. Then, there exists $\delta :=\delta(R,k,s)>0,$ such that for every $x, y \in \mathcal{A}$, we have
\begin{equation}\label{***}
  \big|  |x|^\frac{k}{2}-|y|^\frac{k}{2}  \big|^2 \leq \delta |x-y|^2.
\end{equation}
(In fact, $\delta=2^{k-4} k^2 R^{k-2}$).

\end{lemme}

\begin{preuve}
If $k=2$, we have
\begin{equation}
    \begin{array}{ll}
  \big| |x|-|y| \big|^2 & \leq   |x-y|^2.
    \end{array}
\end{equation}
Let $k>2 $. Let $f(x) :=|x|^\frac{k}{2}$, $x\in\mathbb{R}^n$. We have for $x\in \mathbb{R}^n$, $Df(x)=\frac{k}{2}|x|^{\frac{k}{2}-2}x.$ By the inequality of finite increments applied to the function $f$  we have for $x,y\in\mathbb{R}^n$
\begin{equation}\label{maj x k}
\bigg | |y|^\frac{k}{2}-|x|^\frac{k}{2} \bigg| \leq\frac{k}{2} \displaystyle\sup_{t \in [0,1]} |y+t(x-y)|^{\frac{k}{2}-1} |x-y|,
\end{equation}
Let $H(x,y) :=  \frac{\big|  |x|^\frac{k}{2}-|y|^\frac{k}{2}  \big|^2} {|x-y|^2},$ $x,y\in \mathcal{A}$. Our aim is to prove that $H$ is bounded by some positive constant $\delta$.\\
we obtain for $x,y\in\mathbb{R}^n$
\begin{equation}
\bigg | |y|^\frac{k}{2}-|x|^\frac{k}{2} \bigg|^2 \leq\frac{k^2}{4} \displaystyle\sup_{t \in [0,1]} |y+t(x-y)|^{k-2} |x-y|^2.
\end{equation}
Since for $t \in [0,1]$, $x,y\in \mathcal{A}$, we have  $|y+t(x-y)|\leq 2R$, then
$$\bigg | |y|^\frac{k}{2}-|x|^\frac{k}{2} \bigg|^2 \leq\frac{k^2}{4} (2R)^{k-2} |x-y|^2.$$
Using the fact that $k > 2$, we have for $x,y \in \mathcal{A}$,

\begin{equation}
    \begin{array}{ll}
H(x,y) & \leq\frac{k^2}{4} (2R)^{k-2} \\
\\
&= 2^{k-4} k^2 R^{k-2}
\\
\\
& := \delta.
    \end{array}
\end{equation}
Therefore, $H$ is bounded over $\mathcal{A}$ and the result follows.

\end{preuve}

We move,  now,  to prove \eqref{2} and we take $|x-y|<\gamma<1$. \\ 
We apply the previous lemma for $R=\frac{\eta}{\varepsilon}$

\begin{equation}\label{intégrale sur B}
\begin{array}{ll}
 \displaystyle\int_{|x|<\frac{\eta}{\varepsilon}} \displaystyle\int_{|y|<\frac{\eta}{\varepsilon}}  \frac{\big| |x|^\frac{k}{2}-|y|^\frac{k}{2} \big|^2 1_{B(0,\gamma)}(x-y)}{|x-y|^{n+2s} (1+|y|^2)^{n-2s} } dy dx
 
& \leq \delta \displaystyle\int_{{\mathbb{R}^n} \times{\mathbb{R}^n}} \frac{1}{(1+|y|^2)^{n-2s}}\frac{1_{B(0, \gamma)}(x-y) }{|x-y|^{n+2s-2}}dx dy  \\
\\
& = 2^{k-4} k^2 \eta^{k-2} \varepsilon^{-k+2} \displaystyle\int_{\mathbb{R}^n} \tilde{h}_1*g_1(z) dz  
\end{array}
\end{equation}
where $h_1$ and $g_1$ are the functions defined by 
\begin{equation}
\tilde{h}_1(z):=\frac{1}{(1+|z|^2)^{n-2s}},\  \hbox{and}\  g_1(z) := \frac{1_{B(0, \gamma)}(z)}{|z|^{n+2s-2}}. 
\end{equation}
For $n>4s$, we have ${\tilde{h}}_1 \in L^1(\mathbb{R}^n)$, in other words, provided that 
\begin{equation}\label{contrainte 1}
    \begin{array}{ll}
    \hbox{ if } n=3 \hbox{ and } s<\frac{3}{4}
    \\
    \\
    \hbox{ or } n \geq 4 \hbox{ and } 0<s<1,
    \end{array}
\end{equation} 
which is satisfied since $n$ and $s$ satisfy \eqref{conditions n s}.
Since $s<1$, then we have $g_1\in L^1(\mathbb{R}^n)$. Therefore $\tilde{h}_1*g_1 \in L^1(\mathbb{R}^n)$ and we get
 \begin{equation}\label{premiere partie finie}
   \displaystyle\int_{|x|<\frac{\eta}{\varepsilon}} \displaystyle\int_{|y|<\frac{\eta}{\varepsilon}}  \frac{\big| |x|^\frac{k}{2}-|y|^\frac{k}{2} \big|^2 1_{B(0,\gamma)}(x-y)}{|x-y|^{n+2s} (1+|y|^2)^{n-2s} } dy dx\leq \tilde{C}_0 \varepsilon^{-k+2},
 \end{equation}
 where $\tilde{C}_0$ is a positive constant depending on $k, \eta, n$ and $s.$
 Now, to finish the proof of \eqref{2}, it remains to study the integral in \eqref{2} over the set 
 $\Gamma_\varepsilon:=\{x,y \in \mathbb{R}^n \hbox{ such that } |x|<\frac{\eta}{\varepsilon}, |y|<\frac{\eta}{\varepsilon} \hbox{ and } |x-y|>\gamma \}$.
 Let
\begin{equation}
 X_\varepsilon := \displaystyle\int_{\Gamma_\varepsilon}   \displaystyle\frac{\big| |x|^\frac{k}{2}-|y|^\frac{k}{2} \big|^2}{(1+|y|^2)^{n-2s} |x-y|^{n+2s} } dy dx.
\end{equation}

We have
\begin{equation}\label{x epsilon 1}
\begin{array}{ll}
 X_\varepsilon & = \displaystyle\int_{|x|\leq 1 }\displaystyle\int_{|y|\leq \frac{\eta}{\varepsilon}}  \frac{\big| |x|^\frac{k}{2}-|y|^\frac{k}{2} \big|^2}{(1+|y|^2)^{n-2s} |x-y|^{n+2s} } 1_{B(0,\gamma)^c}(x-y) dy dx \\
 \\
 & +   \displaystyle\int_{1\leq|x|\leq\frac{\eta}{\varepsilon}} \displaystyle\int_{|y|\leq \frac{\eta}{\varepsilon}}   \displaystyle\frac{\big| |x|^\frac{k}{2}-|y|^\frac{k}{2} \big|^2}{(1+|y|^2)^{n-2s} |x-y|^{n+2s} } 1_{B(0,\gamma)^c}(x-y)  dy dx \\
 \end{array}
\end{equation}

 term on the right hand side of \eqref{x epsilon 1}. Since $\displaystyle\frac{1}{|x-y|^{n+2s}} \leq \displaystyle\frac{1}{\gamma^{n+2s}}$ and $k<n-4s$, we have
 
\begin{equation}\label{x inf à 1}
\begin{array}{ll}
     \displaystyle\int_{|x|\leq 1 }\displaystyle\int_{|y|\leq \frac{\eta}{\varepsilon}} 
      \displaystyle\frac{\big| |x|^\frac{k}{2}-|y|^\frac{k}{2} \big|^2}{(1+|y|^2)^{n-2s} |x-y|^{n+2s} } 1_{B(0,\gamma)^c}(x-y) dy dx \\
     \\
      \leq C  \displaystyle\int_{|x|\leq 1 }dx\displaystyle\int_{|y|\leq \frac{\eta}{\varepsilon}} 
     \displaystyle\frac{\big( 1+|y|^\frac{k}{2} \big)^2}{(1+|y|^2)^{n-2s} } dy 
     \\
     \\
     \leq C'_1, 
     \end{array}
\end{equation}
where $C$ and $C'_1$ are positive constants depending on $n, s, k$ and $\gamma$.
\\
For the second term on the right hand side of \eqref{x epsilon 1} we have

\begin{equation}\label{somme_0}
\begin{array}{ll}
    \displaystyle\int_{1\leq|x|\leq\frac{\eta}{\varepsilon}} \displaystyle\int_{|y|\leq \frac{\eta}{\varepsilon}}   \displaystyle\frac{\big| |x|^\frac{k}{2}-|y|^\frac{k}{2} \big|^2}{(1+|y|^2)^{n-2s} |x-y|^{n+2s} } 1_{B(0,\gamma)^c}(x-y)  dy dx
    \\
    \\
    \leq 
    \displaystyle\int_{1\leq|x|\leq\frac{\eta}{\varepsilon}} \displaystyle\int_{0\leq|y|\leq 1}   \displaystyle\frac{\big| |x|^\frac{k}{2}-|y|^\frac{k}{2} \big|^2}{(1+|y|^2)^{n-2s} |x-y|^{n+2s} } 1_{B(0,\gamma)^c}(x-y)  dy dx
    \\
    \\
    + \displaystyle\int_{1\leq|x|\leq\frac{\eta}{\varepsilon}} \displaystyle\int_{1<|y|\leq \frac{\eta}{\varepsilon}}   \displaystyle\frac{\big| |x|^\frac{k}{2}-|y|^\frac{k}{2} \big|^2}{(1+|y|^2)^{n-2s} |x-y|^{n+2s} } 1_{B(0,\gamma)^c}(x-y)  dy dx.
    \end{array}
\end{equation}
We treat both of the terms of the sum \eqref{somme_0} and we start by the first term. Since for $x,y \in \mathbb{R}^n$ such that $|y| \leq 1$, we have $|x-y|^{n+2s} \geq |x| - |y| \geq (|x| - 1)^{n+2s}$, then

\begin{equation}\label{intégrales 1}
    \begin{array}{ll}
      \displaystyle\int_{1\leq|x|\leq\frac{\eta}{\varepsilon}} \displaystyle\int_{0\leq|y|\leq 1}   \displaystyle\frac{\big| |x|^\frac{k}{2}-|y|^\frac{k}{2} \big|^2}{(1+|y|^2)^{n-2s} |x-y|^{n+2s} } 1_{B(0,\gamma)^c}(x-y)  dy dx
      \\
      \\
      = \displaystyle\int_{1\leq|x|\leq 5 } \displaystyle\int_{0\leq|y|\leq 1}   \displaystyle\frac{\big| |x|^\frac{k}{2}-|y|^\frac{k}{2} \big|^2}{(1+|y|^2)^{n-2s} |x-y|^{n+2s} } 1_{B(0,\gamma)^c}(x-y)  dy dx
      \\
      \\
      + \displaystyle\int_{5\leq|x|\leq\frac{\eta}{\varepsilon}} \displaystyle\int_{0\leq|y|\leq 1}   \displaystyle\frac{\big| |x|^\frac{k}{2}-|y|^\frac{k}{2} \big|^2}{(1+|y|^2)^{n-2s} |x-y|^{n+2s} } 1_{B(0,\gamma)^c}(x-y)  dy dx
      \\
      \\
    \leq \displaystyle\frac{1}{\gamma^{n+2s}} \displaystyle\int_{1\leq|x|\leq 5} \displaystyle\int_{0\leq|y|\leq 1}   \displaystyle\frac{\big| |x|^\frac{k}{2}-|y|^\frac{k}{2} \big|^2}{(1+|y|^2)^{n-2s}}  dy dx
    \\
    \\
    +  \displaystyle\int_{5\leq|x|\leq\frac{\eta}{\varepsilon}} \displaystyle\int_{0\leq|y|\leq 1}   \displaystyle\frac{\big| |x|^\frac{k}{2}-|y|^\frac{k}{2} \big|^2}{(1+|y|^2)^{n-2s} (|x| - 1)^{n+2s}}   dy dx
    \end{array}
\end{equation}
For the first term of the sum in \eqref{intégrales 1} we have
\begin{equation}
\begin{array}{ll}
    \displaystyle\int_{1\leq|x|\leq 5} \displaystyle\int_{0\leq|y|\leq 1}   \displaystyle\frac{\big| |x|^\frac{k}{2}-|y|^\frac{k}{2} \big|^2}{(1+|y|^2)^{n-2s}}  dy dx 
   &
    \leq 2  \displaystyle\int_{1\leq|x|\leq 5} \displaystyle\int_{0\leq|y|\leq 1}   \displaystyle\frac{ (|x|^k+|y|^k) }{(1+|y|^2)^{n-2s}}  dy dx
    \\
    \\
    & = 2  \displaystyle\int_{1\leq|x|\leq 5} |x|^k dx 
    \displaystyle\int_{0\leq|y|\leq 1} \displaystyle\frac{dy}{(1+|y|^2)^{n-2s}} \\
    \\
    & + 2  \displaystyle\int_{1\leq|x|\leq 5}  dx 
    \displaystyle\int_{0\leq|y|\leq 1} \displaystyle\frac{|y|^k}{(1+|y|^2)^{n-2s}}dy,
    \end{array}
\end{equation}
every integral here represents an integral of a continuous function over a compact set, then each integral is convergent (and is bounded independently of $\varepsilon$) and we get 
\begin{equation}\label{100}
      \displaystyle\int_{1\leq|x|\leq 5} \displaystyle\int_{0\leq|y|\leq 1}   \displaystyle\frac{\big| |x|^\frac{k}{2}-|y|^\frac{k}{2} \big|^2}{(1+|y|^2)^{n-2s}}  dy dx  \leq C,
\end{equation}
where $C$ is a positive constant depending on $n,s$ and $k.$\\
Now, for the second term of the sum in \eqref{intégrales 1}, we write
\begin{equation}\label{5}
    \begin{array}{ll}
    \displaystyle\int_{5\leq|x|\leq\frac{\eta}{\varepsilon}} \displaystyle\int_{0\leq|y|\leq 1}   \displaystyle\frac{\big| |x|^\frac{k}{2}-|y|^\frac{k}{2} \big|^2}{(1+|y|^2)^{n-2s} (|x| - 1)^{n+2s}}   dy dx
    \\
    \\
    \leq 2  \displaystyle\int_{5\leq|x|\leq\frac{\eta}{\varepsilon}} \displaystyle\int_{0\leq|y|\leq 1} 
       \displaystyle\frac{( |x|^k+|y|^k) }{(1+|y|^2)^{n-2s} (|x| - 1)^{n+2s}}   dy dx
       \\
       \\
       \leq 2   \displaystyle\int_{5\leq|x|\leq\frac{\eta}{\varepsilon}} \displaystyle\frac{|x|^k}{(|x|-1)^{n+2s}} dx 
       \displaystyle\int_{0\leq|y|\leq 1} \displaystyle\frac{dy}{(1+|y|^2)^{n-2s} } 
       \\
       \\
       + 2  \displaystyle\int_{5\leq|x|\leq +\infty} 
       \displaystyle\frac{dx}{(|x|-1)^{n+2s}} 
       \displaystyle\int_{0\leq|y|\leq 1}  \displaystyle\frac{|y|^k}{(1+|y|^2)^{n-2s}}dy.
    \end{array}
\end{equation}
Clearly, the integrals $ \displaystyle\int_{0\leq|y|\leq 1} \displaystyle\frac{dy}{(1+|y|^2)^{n-2s} } $ and   $\displaystyle\int_{0\leq|y|\leq 1}  \displaystyle\frac{|y|^k}{(1+|y|^2)^{n-2s}}dy$ are convergent (they are bounded independently of $\varepsilon$). \\
Since $2s<k$, we have 
\begin{equation}\label{3}
 \displaystyle\int_{5\leq|x|\leq\frac{\eta}{\varepsilon}} \displaystyle\frac{|x|^k}{(|x|-1)^{n+2s}} dx \leq \bar{C} \varepsilon^{-k+2s},
\end{equation}
where $\bar{C}$ is a positive constant possibly depending on $n,s$ and $k$
Since $s>0$, we get
\begin{equation}\label{4}
    \displaystyle\int_{5\leq|x|\leq +\infty} 
       \displaystyle\frac{dx}{(|x|-1)^{n+2s}} \leq C,
\end{equation}
where $C$ is a positive constant depending only on $n$ and $s$.\\
Therefore, using \eqref{3} and \eqref{4} together with \eqref{5} we get for $2s<k<n-4s,$ $s\in ]0,1[$ and $\varepsilon$ sufficiently small
\begin{equation}\label{epsilon_10} 
 \displaystyle\int_{5\leq|x|\leq\frac{\eta}{\varepsilon}} \displaystyle\int_{0\leq|y|\leq 1}   \displaystyle\frac{\big| |x|^\frac{k}{2}-|y|^\frac{k}{2} \big|^2}{(1+|y|^2)^{n-2s} (|x| - 1)^{n+2s}}   dy dx \leq  C  \varepsilon^{-k+2s},
 \end{equation}
where $C$ is a positive constant possibly depending on $n,s$ and $k.$
\\
Thus, thanks to \eqref{epsilon_10}, \eqref{100} together with \eqref{intégrales 1} we have

\begin{equation}\label{intégrale y entre 0 et 1}
\begin{array}{ll}
     \displaystyle\int_{1\leq|x|\leq\frac{\eta}{\varepsilon}} \displaystyle\int_{0\leq|y|\leq 1}   \displaystyle\frac{\big| |x|^\frac{k}{2}-|y|^\frac{k}{2} \big|^2}{(1+|y|^2)^{n-2s} |x-y|^{n+2s} } 1_{B(0,\gamma)^c}(x-y)  dy dx
     & 
     \leq C \varepsilon^{-k+2s}
     \end{array}
\end{equation}
where $ C$ is a positive constant possibly depending on $n,s$ and $k.$\\
Let us consider the second term of the sum \eqref{somme_0}. Using an argument in \cite{FV}, for $x,y\in \mathbb{R}^n$, let $x= rx'$, $y=\rho y'$, and $|x'|=|y'|=1,$ and perform the following change of variable: $\tau=\displaystyle\frac{\rho}{r}$ we have 
\begin{equation}\label{intégrale 2}
    \begin{array}{ll}
       \displaystyle\int_{1\leq|x|\leq\frac{\eta}{\varepsilon}} \displaystyle\int_{1<|y|\leq \frac{\eta}{\varepsilon}}   \displaystyle\frac{\big| |x|^\frac{k}{2}-|y|^\frac{k}{2} \big|^2}{(1+|y|^2)^{n-2s} |x-y|^{n+2s} } 1_{B(0,\gamma)^c}(x-y)  dy dx
       \\
       \\
       \leq \displaystyle\int_1^{\frac{\eta}{\varepsilon}} \displaystyle\int_1^\frac{\eta}{\varepsilon} \displaystyle\int_{|y'|=1} \displaystyle\int_{|x'|=1} \displaystyle\frac{|\rho^\frac{k}{2}-r^\frac{k}{2}|^2}{(1+\rho^2)^{n-2s}|rx'-\rho y'|^{n+2s}} d\sigma(x') d\sigma(y') \rho^{n-1} d\rho r^{n-1} dr
       \\
       \\
       \leq  \displaystyle\int_1^{\frac{\eta}{\varepsilon}} \displaystyle\int_1^\frac{\eta}{\varepsilon} \displaystyle\int_{|y'|=1} \displaystyle\int_{|x'|=1} \displaystyle\frac{r^k |1-(\frac{\rho}{r})^\frac{k}{2}|^2 r^{n-1}\rho^{n-1}}{\rho^{2n-4s} r^{n+2s} |x'-\frac{\rho}{r}y'|^{n+2s}} d\sigma(x') d\sigma(y') d\rho  dr
       \\
       \\
      = \displaystyle\int_1^{\frac{\eta}{\varepsilon}}       r^{k-2s-1}\displaystyle\int_\frac{1}{r}^\frac{\eta}{\varepsilon r} \displaystyle\int_{|y'|=1} \displaystyle\int_{|x'|=1} \displaystyle\frac{|1-\tau^\frac{k}{2}|^2 \tau^{n-1}r^{n-1}r}{r^{2n-4s} \tau^{2n-4s}|x'-\tau y'|^{n+2s}} 
       d\sigma(x') d\sigma(y') d\tau  dr
       \\
       \\
    =   \displaystyle\int_1^{\frac{\eta}{\varepsilon}}
       r^{k+2s-n-1}\displaystyle\int_\frac{1}{r}^\frac{1}{2} \displaystyle\int_{|y'|=1} \displaystyle\int_{|x'|=1} 
       \displaystyle\frac{|1-\tau^\frac{k}{2}|^2 \tau^{-n+4s-1}}{|x'-\tau y'|^{n+2s}}
        d\sigma(x') d\sigma(y') d\tau  dr
        \\
        \\
        +  \displaystyle\int_1^{\frac{\eta}{\varepsilon}}
       r^{k+2s-n-1}\displaystyle\int_\frac{1}{2}^1 \displaystyle\int_{|y'|=1} \displaystyle\int_{|x'|=1} 
       \displaystyle\frac{|1-\tau^\frac{k}{2}|^2 \tau^{-n+4s-1}}{|x'-\tau y'|^{n+2s}}
        d\sigma(x') d\sigma(y') d\tau  dr
        \\
        \\
 + \displaystyle\int_1^{\frac{\eta}{\varepsilon}}
       r^{k+2s-n-1}\displaystyle\int_1^\frac{\eta}{\varepsilon r} \displaystyle\int_{|y'|=1} \displaystyle\int_{|x'|=1} 
       \displaystyle\frac{|1-\tau^\frac{k}{2}|^2 \tau^{-n+4s-1}}{|x'-\tau y'|^{n+2s}}
        d\sigma(x') d\sigma(y') d\tau  dr.
     \end{array}
\end{equation}
We will estimate the terms of the previous sum and we begin by the first term. Since $|x-y|>\gamma,$ then there exists $\bar{\delta}>0$ such that $|x'-\tau y'|> \bar{\delta}$ and since $2s<k<n-2s$ (in fact, $k \in [2,n-4s[$), for $\varepsilon$ small enough we have 
\begin{equation}\label{first term}
    \begin{array}{ll}
   \displaystyle\int_1^{\frac{\eta}{\varepsilon}}
       r^{k+2s-n-1}\displaystyle\int_\frac{1}{r}^\frac{1}{2} \displaystyle\int_{|y'|=1} \displaystyle\int_{|x'|=1} 
       \displaystyle\frac{|1-\tau^\frac{k}{2}|^2 \tau^{-n+4s-1}}{|x'-\tau y'|^{n+2s}}
        d\sigma(x') d\sigma(y') d\tau  dr
        \\
        \\
        \leq C \displaystyle\int_1^{\frac{\eta}{\varepsilon}}
       r^{k+2s-n-1}\displaystyle\int_\frac{1}{r}^\frac{1}{2}  
       |1-\tau^\frac{k}{2}|^2 \tau^{-n+4s-1} d\tau  dr
       \\
       \\
       \leq C \displaystyle\int_1^{\frac{\eta}{\varepsilon}}
       r^{k+2s-n-1}\displaystyle\int_\frac{1}{r}^\frac{1}{2} (\tau^{-n+4s-1}+\tau^{-n+4s-1+k}) d\tau dr
       \\
       \\
       \leq C \bigg[ (\eta^{k-2s} \varepsilon^{-k+2s}-1)
    +(1-\eta^{k+2s-n} \varepsilon^{-k-2s+n}) 
    + (1-\eta^{-2s} \varepsilon^{2s})\bigg]
    \\
    \\
    \leq C \varepsilon^{-k+2s},
    \end{array}
\end{equation}
where $C$ is a positive constant possibly depending on $n, s, k$ and $\eta.$
\\
\\
Notice that $ K(\tau):=\displaystyle\int_{|y'|=1}  \displaystyle\frac{1}{|x'-\tau y'|^{n+2s}}   d\sigma(y')$ is independent of $x' \in \{ |y|=1 \}$ (see \cite{FV}) and we have $K(\tau)= \tau^{2-n} (\tau^2-1)^{-1-2s} H(\tau),$ $\tau \geq1$, where $H(\tau)$ is a positive continuous function on $[1,+\infty[$ such that $H(\tau) \approx \tau^{2s}$ as $\tau$ tends to $+\infty$. 
\\
Therefore for the second term of the sum in \eqref{intégrale 2}, performing a change of variables $\xi := \frac{1}{\tau}$, we have 
\begin{equation}\label{somme trois termes}
\begin{array}{ll}
     \displaystyle\int_1^{\frac{\eta}{\varepsilon}}
       r^{k+2s-n-1}\displaystyle\int_\frac{1}{2}^1 \displaystyle\int_{|y'|=1} \displaystyle\int_{|x'|=1} 
       \displaystyle\frac{|1-\tau^\frac{k}{2}|^2 \tau^{-n+4s-1}}{|x'-\tau y'|^{n+2s}}
        d\sigma(x') d\sigma(y') d\tau  dr
        \\
        \\
        =   \displaystyle\int_1^{\frac{\eta}{\varepsilon}}
       r^{k+2s-n-1}\displaystyle\int_\frac{1}{2}^1\displaystyle\int_{|x'|=1} \bigg(\displaystyle\int_{|y'|=1} \displaystyle\frac{1}{|x'-\tau y'|^{n+2s}}  d\sigma(y')\bigg) 
      |1-\tau^\frac{k}{2}|^2 \tau^{-n+4s-1}
        d\sigma(x') d\tau  dr
        \\
        \\
        = \sigma(S^{n-1}) \displaystyle\int_1^{\frac{\eta}{\varepsilon}}
       r^{k+2s-n-1}\displaystyle\int_\frac{1}{2}^1 K(\tau)
      |1-\tau^\frac{k}{2}|^2 \tau^{-n+4s-1}
         d\tau  dr
        \\
        \\
         =\sigma(S^{n-1}) \displaystyle\int_1^{\frac{\eta}{\varepsilon}}
       r^{k+2s-n-1}dr\displaystyle\int_1^2 |1-\xi^\frac{-k}{2}|^2 \xi^{n-4s-1} K(\frac{1}{\xi}) d\xi
       \\
       \\
       = \sigma(S^{n-1}) \displaystyle\int_1^{\frac{\eta}{\varepsilon}}
       r^{k+2s-n-1}dr\displaystyle\int_1^2 |1-\xi^\frac{-k}{2}|^2 \xi^{n-4s-1} \xi^{n+2s} K(\xi) d\xi
       \\
       \\
       = \sigma(S^{n-1}) \displaystyle\int_1^{\frac{\eta}{\varepsilon}}
       r^{k+2s-n-1}dr\displaystyle\int_1^2 |1-\xi^\frac{-k}{2}|^2 \xi^{2n-2s-1}  H(\xi) \xi^{2-n} (\xi^2-1)^{-1-2s} d\xi
       \\
       \\
       = \sigma(S^{n-1}) \displaystyle\int_1^{\frac{\eta}{\varepsilon}}
       r^{k+2s-n-1}dr\displaystyle\int_1^2 |1-\xi^\frac{-k}{2}|^2 \xi^{n-2s+1}  H(\xi)  (\xi^2-1)^{-1-2s} d\xi,
        \end{array}
\end{equation}
where $ \sigma(S^{n-1})$ is the Lebesgue measure of the unit sphere in $\mathbb{R}^n$.
Since $k<n-4s$, then $k<n-2s$ and we have
\begin{equation}\label{1ère intég*}
    \begin{array}{ll}
    \displaystyle\int_1^\frac{\eta}{\varepsilon }  
       r^{k+2s-n-1} dr

     \leq \displaystyle\frac{1}{n-2s-k} := C_4,
    \end{array}
\end{equation}
clearly $C_4$ is a positive constant depending on $n,s$ and $k.$\\
The function $ |1-\xi^\frac{-k}{2}|^2 \xi^{n-2s+1}  H(\xi)  (\xi^2-1)^{-1-2s} $ is continuous on $]1,2]$ and in the neighborhood of $1$ this functions is equivalent to 
$\displaystyle\frac{|1-\xi^{\frac{-k}{2}}|^2 H(1)}{(\xi^2-1)^{1+2s}} = \xi^{-k}  \displaystyle\frac{|\xi^\frac{k}{2}-1|^2 H(1)}{(\xi^2-1)^{1+2s}}$
which in the neighborhood of $\xi=1$ is equivalent to
$ \displaystyle\frac{|\xi^\frac{k}{2}-1|^2 H(1)}{(\xi^2-1)^{1+2s}}$ and if we note $\xi':=\xi-1,$  we get after performing a limited development on $\xi'=0$
\begin{equation}\label{dev limit}
\begin{array}{ll}
   \displaystyle\frac{(\xi^\frac{k}{2}-1)^2}{(\xi^2-1)^{1+2s}} 
   & = \displaystyle\frac{\big(  (\xi'+1)^\frac{k}{2}-1  \big)^2}{\big( (\xi'+1)^2-1  \big)^{1+2s}} 
   \\
   \\
   & = \xi'^{(1-2s)} \bigg(\displaystyle\frac{(\frac{k}{2})^2}{2^{1+2s}}+o(1)\bigg),
\end{array}
\end{equation}
which is integrable in the neighborhood of $\xi'=0$ since $s<1.$\\
Therefore we have
\begin{equation}\label{C_5}
    \displaystyle\int_1^2 |1-\xi^\frac{-k}{2}|^2 \xi^{n-2s+1}  H(\xi)  (\xi^2-1)^{-1-2s} d\xi \leq C_5,
\end{equation}
where $C_5$ is a positive constant depending on $s,k$ and $n.$\\
Thus, using \eqref{C_5} and \eqref{1ère intég*} together with \eqref{somme trois termes} we find
\begin{equation}\label{second term}
      \displaystyle\int_1^{\frac{\eta}{\varepsilon}}
       r^{k+2s-n-1}\displaystyle\int_\frac{1}{2}^1 \displaystyle\int_{|y'|=1} \displaystyle\int_{|x'|=1} 
       \displaystyle\frac{|1-\tau^\frac{k}{2}|^2 \tau^{-n+4s-1}}{|x'-\tau y'|^{n+2s}}
        d\sigma(x') d\sigma(y') d\tau  dr \leq C_6,
\end{equation}
where $C_6$ is a positive constant depending on $s,k$ and $n.$ 
\\
We treat the third term of the sum in \eqref{intégrale 2} and we proceed as in \eqref{somme trois termes} we have
\begin{equation}\label{K}
    \begin{array}{ll}
       \displaystyle\int_1^{\frac{\eta}{\varepsilon}}
       r^{k+2s-n-1}\displaystyle\int_1^\frac{\eta}{\varepsilon r} \displaystyle\int_{|y'|=1} \displaystyle\int_{|x'|=1} 
       \displaystyle\frac{|1-\tau^\frac{k}{2}|^2 \tau^{-n+4s-1}}{|x'-\tau y'|^{n+2s}}
        d\sigma(x') d\sigma(y') d\tau  dr
       \\
       \\
      
     =  \displaystyle\int_1^{\frac{\eta}{\varepsilon}}
       r^{k+2s-n-1}\displaystyle\int_1^\frac{\eta}{\varepsilon r} \displaystyle\int_{|x'|=1}  K(\tau)
      |1-\tau^\frac{k}{2}|^2 \tau^{-n+4s-1} 
        d\sigma(x') d\tau  dr
        \\
        \\
     = \sigma(S^{n-1})\displaystyle\int_1^{\frac{\eta}{\varepsilon}}
       r^{k+2s-n-1}\displaystyle\int_1^\frac{\eta}{\varepsilon r}   H(\tau) \tau^{2-n} (\tau^2-1)^{-1-2s} 
      |1-\tau^\frac{k}{2}|^2 \tau^{-n+4s-1} 
        d\tau  dr
        \\
        \\
       \leq  \sigma(S^{n-1}) \displaystyle\int_1^\frac{\eta}{\varepsilon }  
       r^{k+2s-n-1} dr \displaystyle\int_1^{+\infty}\tau^{-2n+4s+1} H(\tau)  (\tau^2-1)^{-1-2s} 
      |1-\tau^\frac{k}{2}|^2
       d\tau  
    \end{array}
\end{equation}
where $ \sigma(S^{n-1})$ is the Lebesgue measure of the unit sphere in $\mathbb{R}^n$.\\
Let us treat both of the integrals in the last expression of \eqref{K}, for $k<n-2s$ (which is satisfied since $k<n-4s$), by \eqref{1ère intég*} we have
\begin{equation}\label{1ère intég}
    \begin{array}{ll}
    \displaystyle\int_1^\frac{\eta}{\varepsilon }  
       r^{k+2s-n-1} dr \leq C_4.
    \end{array}
\end{equation}
\\
The second integral of \eqref{K} which is the following
\begin{equation}\label{intégrale11}
   \displaystyle\int_1^{+\infty}\tau^{-2n+4s+1} H(\tau)  (\tau^2-1)^{-1-2s} 
      |1-\tau^\frac{k}{2}|^2
       d\tau ,
\end{equation}
it is the integral of a continuous function on  $]1, +\infty[$. In the neighborhood of $+\infty$, since  $H(\tau) \approx \tau^{2s}$ as $\tau$ tends to $+\infty$, the function inside the integral is equivalent to 
 $\displaystyle\frac{1}{\tau^{2n-k-2s+1}},$
which is integrable at $+\infty$ since $k<2n-2s$ (in fact, $k<n-4s$ and $0<s<1$).\\
In the neighborhood of $1,$ after performing the changement $\tau'= \tau-1$, the function inside the integral \eqref{intégrale11} is equivalent to
$   \displaystyle\frac{(\tau^\frac{k}{2}-1)^2}{(\tau^2-1)^{1+2s}} ,$ we use the same computations as in \eqref{dev limit} and since $s<1$ we find that it
is integrable in the neighborhood of $\tau=1$. 
\\
Then, we get 
\begin{equation}\label{2eme intég}
      \displaystyle\int_1^{+\infty}\tau^{-2n+4s+1} H(\tau)  (\tau^2-1)^{-1-2s} 
      |1-\tau^\frac{k}{2}|^2
       d\tau \leq C_1,
\end{equation}
where $C_1$ is a positive constant depending on $n,s$ and $k$.
Thus, using \eqref{1ère intég} and \eqref{2eme intég} together with 
\eqref{K} we have 
\begin{equation}\label{third term}
\begin{array}{ll}
       \displaystyle\int_1^{\frac{\eta}{\varepsilon}}
       r^{k+2s-n-1}\displaystyle\int_1^\frac{\eta}{\varepsilon r} \displaystyle\int_{|y'|=1} \displaystyle\int_{|x'|=1} 
       \displaystyle\frac{|1-\tau^\frac{k}{2}|^2 \tau^{-n+4s-1}}{|x'-\tau y'|^{n+2s}}
        d\sigma(x') d\sigma(y') d\tau  dr
      & \leq \tilde{C}_2,
\end{array}
\end{equation}
where $\tilde{C_2}$ is a positive constant depending on $n,s$ and $k$. 
Thus, thanks to \eqref{first term}, \eqref{second term}, \eqref{third term} together with \eqref{intégrale 2} we get
\begin{equation} \label{intégrales}
\begin{array}{ll}
      \displaystyle\int_{1\leq|x|\leq\frac{\eta}{\varepsilon}} \displaystyle\int_{1<|y|\leq \frac{\eta}{\varepsilon}}   \displaystyle\frac{\big| |x|^\frac{k}{2}-|y|^\frac{k}{2} \big|^2}{(1+|y|^2)^{n-2s} |x-y|^{n+2s} } 1_{B(0,\gamma)^c}(x-y)  dy dx 
\leq C'_1+C'_2 \varepsilon^{-k+2s},
\end{array}
\end{equation}
where $C'_1$ and $C'_2$ are positive constants possibly depending on $n,s,\eta$ and $k.$
\\
Then, by \eqref{intégrale y entre 0 et 1} and \eqref{intégrales} together with \eqref{somme_0} we get
\begin{equation}\label{x supérieure à 1}
\begin{array}{ll}

      \displaystyle\int_{1\leq|x|\leq\frac{\eta}{\varepsilon}} \displaystyle\int_{|y|\leq \frac{\eta}{\varepsilon}}   \displaystyle\frac{\big| |x|^\frac{k}{2}-|y|^\frac{k}{2} \big|^2}{(1+|y|^2)^{n-2s} |x-y|^{n+2s} } 1_{B(0,\gamma)^c}(x-y)  dy dx  \leq C''_1 \varepsilon^{-k+2s} + C_3,
      \end{array}
\end{equation}
where $C''_1$ and $C_3$ are positive constants possibly depending on $n,s,\eta$ and $k.$
\\
Therefore, using \eqref{x inf à 1} and \eqref{x supérieure à 1}  together with \eqref{x epsilon 1}, we get for $2s<k<n-4s$
\begin{equation}\label{deuxième partie finie}
\begin{array}{ll}
     X_\varepsilon & = \displaystyle\int_{\Gamma_\varepsilon}   \displaystyle\frac{\big| |x|^\frac{k}{2}-|y|^\frac{k}{2} \big|^2}{(1+|y|^2)^{n-2s} |x-y|^{n+2s} } dy dx
     \leq \tilde{C}_1+ \tilde{C}_3  \varepsilon^{-k+2s},
     \end{array}
\end{equation}
where $\tilde{C}_1$ and $\tilde{C}_3$ are  positive constants possibly depending on $n,s, \eta$ and $k.$ \\
Thus, by \eqref{premiere partie finie} and \eqref{deuxième partie finie}, for $2s<k<n-4s$
\begin{equation}
    \begin{array}{ll}
\displaystyle\int_{|x|\leq\frac{\eta}{\varepsilon}} \displaystyle\int_{|y|\leq\frac{\eta}{\varepsilon}}\frac{\big| |x|^\frac{k}{2}-|y|^\frac{k}{2} \big|^2}{(1+|y|^2)^{n-2s} |x-y|^{n+2s}} dy dx \\
\\ = \displaystyle\int_{|x|\leq\frac{\eta}{\varepsilon}} \displaystyle\int_{|y|\leq\frac{\eta}{\varepsilon}}\frac{\big| |x|^\frac{k}{2}-|y|^\frac{k}{2} \big|^2}{(1+|y|^2)^{n-2s} |x-y|^{n+2s}  } 1_{B(0,\gamma)}(x-y) dy dx \\
\\
+ \displaystyle\int_{\Gamma_\varepsilon}  \frac{\big| |x|^\frac{k}{2}-|y|^\frac{k}{2} \big|^2}{(1+|y|^2)^{n-2s} |x-y|^{n+2s} } dy dx\\
\\
 \leq C'+ C" \varepsilon^{-k+2s}+\tilde{C}_0 \varepsilon^{-k+2},
      \end{array} 
\end{equation}
where $C'$ and $C"$ are positive constants possibly depending on $n,s, \eta$ and $k.$ Thus \eqref{2} yields.\\
By \eqref{inégalité clés}, \eqref{1} and \eqref{2}, we finally conclude the following assertion:\\
For all $n$ and $s$ as described in \eqref{conditions n s} and for  $2 \leq k<n-4s,$ with taking $\varepsilon$ small enough we obtain 
\begin{equation}
    \begin{array}{ll}
    A_{s, k, \varepsilon} & =\displaystyle\int_{|x|\leq \eta} |x|^k \bigg ( \displaystyle\int_{|y|\leq \eta}\frac{|u_{\varepsilon,s}(x)-u_{\varepsilon,s}(y)|^2}{|x-y|^{n+2s}}dy\Bigg ) dx \\
    \\ 
    & \leq 2 \varepsilon^k [\bar{\tilde{C}} + (C'+ C" \varepsilon^{-k+2s}+\tilde{C}_0 \varepsilon^{-k+2})]
    \\
    \\
    & \leq C \varepsilon^{2s},
    \end{array}
\end{equation}
where $C$ is a positive constant depending on $s,k,\gamma$ and $n.$
\end{preuve}\\
We move,  now,  to the fourth term in the right hand side of \eqref{4 intégrales}. 
Using the assertion \eqref{b)},  we have
\begin{equation}
    \begin{array}{ll}
       \displaystyle\int_{B(0, \eta)^c}\displaystyle\int_ {B(0, \eta)^c} (p(x)-p_0) \frac{|u_{\varepsilon,s}(x)-u_{\varepsilon,s}(y)|^2}{|x-y|^{n+2s}}dx dy & \\ 
       \\
       \leq C \varepsilon^{n-2s}\displaystyle\int_{B(0, \eta)^c} \bigg( \displaystyle\int_{B(0, \eta)^c}(p(x)-p_0)\frac{\min\{1, |x-y|^{2}\}}{|x-y|^{n+2s}} dx \bigg) dy &\\
       \\
       = C \varepsilon^{n-2s}\bigg[ \displaystyle\int_{B(0, \eta)^c}\bigg(\displaystyle\int_{B(0, \eta)^c} (p(x)-p_0)\frac{1_{B(0,1)}(x-y)}{|x-y|^{n+2s-2}}dx \bigg) dy 
       &\\
       \\
    +\displaystyle\int_{B(0, \eta)^c}\bigg(\displaystyle\int_{B(0, \eta)^c} (p(x)-p_0) \frac{1_{B(0,1)^c}(x-y)}{|x-y|^{n+2s}}dx \bigg) dy  \bigg] & \\
    \\
    \leq C \varepsilon^{n-2s}\bigg[ \displaystyle\int_{\mathbb{R}^n}\bigg(\displaystyle\int_{B(0, \eta)^c} (p(x)-p_0) \frac{1_{B(0,1)}(x-y)}{|x-y|^{n+2s-2}}dx \bigg) dy 
       &\\
       \\
    +\displaystyle\int_{\mathbb{R}^n}\bigg(\displaystyle\int_{B(0, \eta)^c} (p(x)-p_0) \frac{1_{B(0,1)^c}(x-y)}{|x-y|^{n+2s}}dx \bigg) dy  \bigg] 
    &\\
       \\
        = C \varepsilon^{n-2s} \displaystyle\int_{R^n} f_{1, k}*f_2(y) dy +C \varepsilon^{n-2s} \displaystyle\int_{R^n} f_{1, k}*f_{3}(y) dy,  &\\
    \end{array}
\end{equation}
where for all $z\in \mathbb{R}^n$,  
\begin{equation}
f_{1, k}(z) := (p(z)-p_0) 1_{B(0, \eta)^c}(z), \  f_2(z) :=\displaystyle\frac{1_{B(0, 1)}(z)}{|z|^{n+2s-2}} \ \mbox{and} \ f_3(z):=\displaystyle\frac{1_{B(0, 1)^c}(z)}{|z|^{n+2s}}.
\end{equation}
Using \eqref{p-p_0 dans L1} we have $f_{1, k} \in L^1(\mathbb{R}^n)$.  Since $0<s<1$, the functions $f_2$ and $f_3$ are in $L^1(\mathbb{R}^n)$ ,  then $f_{1, k}*f_2$ and $f_{1, k}*f_3$ are in $L^1(\mathbb{R}^n)$.\\
Therefore
 \begin{equation}\label{45}
 \displaystyle\int_{B(0, \eta)^c}\displaystyle\int_ {B(0, \eta)^c} (p(x)-p_0) \frac{|u_{\varepsilon,s}(x)-u_{\varepsilon,s}(y)|^2}{|x-y|^{n+2s}}dx dy = \mathcal{O}( \varepsilon^{n-2s}).
 \end{equation}
 We move,  now,  to study the integral over the set $\mathbb{E}$.
 By using the assertion \eqref{a)},  we have
 \begin{equation}\label{40}
 \begin{array}{ll}
 \displaystyle\int_ {\mathbb{E}} |x|^k\frac{|u_{\varepsilon,s}(x)-u_{\varepsilon,s}(y)|^2}{|x-y|^{n+2s}}dx dy &\\
 \\
 \leq C \varepsilon^{n-2s} \displaystyle\int_{x\in B(0, \eta), y\in B(0, \eta)^c,  |x-y|\leq\frac{\eta}{2}} |x|^k\frac{|x-y|^2}{|x-y|^{n+2s}}dx dy&\\
 \\
 \leq C \varepsilon^{n-2s} \eta^k \displaystyle\int_{x\in B(0, \eta), y\in B(0, \eta)^c,  |x-y|\leq\frac{\eta}{2}}\frac{1}{|x-y|^{n+2s-2}}dx dy &\\
 \\
  \leq C \varepsilon^{n-2s} \eta^k \displaystyle\int_{|x|<\eta} dx \displaystyle\int_{|\xi|\leq\frac{\eta}{2}}\frac{1}{|\xi|^{n+2s-2}} d\xi &\\
  \\
 = \mathcal{O}(\varepsilon^{n-2s}), &\\
\end{array}
 \end{equation}
 as $\varepsilon$ tends to zero. In both these estimates,  we use that $s\in ]0, 1[$.\\
 Now,  in \eqref{4 intégrales},  it remains to estimate the integral on $\mathbb{D}$,  that is,  
 \begin{equation}\label{D}
\displaystyle\int_{\mathbb{D}} |x|^k\frac{|u_{\varepsilon,s}(x)-u_{\varepsilon,s}(y)|^2}{|x-y|^{n+2s}}dx dy.
\end{equation}
For this, we use the same computation as done for (14,27) in \cite{BRS}. Let us recall that: $u_{\varepsilon,s}(x)=U_{\varepsilon, s}(x), $ for any $x\in B(0, \eta)$.\\
We note that for any $(x, y)\in\mathbb{D}$, 

\begin{equation}
\begin{array}{ll}
 |u_{\varepsilon,s}(x)-u_{\varepsilon,s}(y)|^2 & = |U_{\varepsilon, s}(x)-u_{\varepsilon,s}(y)|^2\\ \\
& =\displaystyle|(U_{\varepsilon, s}(x)-U_{\varepsilon, s}(y))+(U_{\varepsilon, s}(y)-u_{\varepsilon,s}(y))|^2\\\\
& \leq |(U_{\varepsilon, s}(x)-U_{\varepsilon, s}(y))|^2+|U_{\varepsilon, s}(y)-u_{\varepsilon,s}(y)|^2\\ \\
& + 2 |(U_{\varepsilon, s}(x)-U_{\varepsilon, s}(y))| |U_{\varepsilon, s}(y)-u_{\varepsilon,s}(y)|, 
\end{array}
\end{equation}
so
\begin{equation}\label{somme D}
    \begin{array}{ll}
    \displaystyle\int_{\mathbb{D}} |x|^k\frac{|u_{\varepsilon,s}(x)-u_{\varepsilon,s}(y)|^2}{|x-y|^{n+2s}}dx dy & \\
    \\
    \leq \displaystyle\int_{\mathbb{D}} |x|^k\frac{|(U_{\varepsilon, s}(x)-U_{\varepsilon, s}(y))|^2}{|x-y|^{n+2s}}dx dy &\\
    \\
    + \displaystyle\int_{\mathbb{D}} |x|^k\frac{|U_{\varepsilon, s}(y)-u_{\varepsilon,s}(y)|^2}{|x-y|^{n+2s}}dx dy& \\
     \\
      +2 \displaystyle\int_{\mathbb{D}} |x|^k\frac{|U_{\varepsilon, s}(x)-U_{\varepsilon, s}(y)| |U_{\varepsilon, s}(y)-u_{\varepsilon,s}(y)|}{|x-y|^{n+2s}}dx dy.&\\ 
    \end{array}
\end{equation}
Hence,  to estimate \eqref{D},  we bound the three terms on the right hand side of \eqref{somme D}. By exploiting \eqref{majoration fonction test u} (here used with $\rho=\eta$),  we obtain
\begin{equation}\label{2ème D}
\begin{array}{ll}
\displaystyle\int_{\mathbb{D}} |x|^k\frac{|U_{\varepsilon, s}(y)-u_{\varepsilon,s}(y)|^2}{|x-y|^{n+2s}}dx dy &

\leq \displaystyle\int_{\mathbb{D}} |x|^k\frac{(|U_{\varepsilon, s}(y)|+|u_{\varepsilon,s}(y)|)^2}{|x-y|^{n+2s}}dx dy \\
\\
&
  \leq 4 \displaystyle\int_{\mathbb{D}} |x|^k\frac{|U_{\varepsilon, s}(y)|^2}{|x-y|^{n+2s}}dx dy \\
  \\
  &
  \leq C \varepsilon^{n-2s} \displaystyle\int_{x\in B(0, \eta), y\in B(0, \eta)^c,  |x-y|>\frac{\eta}{2}}\frac{|x|^k}{|x-y|^{n+2s}}dx dy\\
  \\
  &
  \leq C \eta^k \varepsilon^{n-2s} \displaystyle\int_{x\in B(0, \eta), y\in B(0, \eta)^c,  |x-y|>\frac{\eta}{2}}\frac{1}{|x-y|^{n+2s}}dx dy\\
  \\
  &
  \leq C \eta^k \varepsilon^{n-2s} \displaystyle\int_{|\zeta| < \eta} d\zeta \displaystyle\int_{|\xi|>\frac{\eta}{2}}\frac{1}{|\xi|^{n+2s}} d\xi\\
   \\
  &
  = \mathcal{O}(\varepsilon^{n-2s}), 
\end{array}
\end{equation}
as $\varepsilon$ tends to zero.

We estimate now the last term on the right-hand side of \eqref{somme D}.

\begin{equation}\label{23}
\begin{array}{ll}
    \displaystyle\int_{\mathbb{D}} |x|^k\frac{|U_{\varepsilon, s}(x)-U_{\varepsilon, s}(y)| |U_{\varepsilon, s}(y)-u_{\varepsilon,s}(y)|}{|x-y|^{n+2s}}dx dy & \\
    \\
    \leq \displaystyle\int_{\mathbb{D}} |x|^k\frac{|U_{\varepsilon, s}(x)| |U_{\varepsilon, s}(y)-u_{\varepsilon,s}(y)|}{|x-y|^{n+2s}}dx dy 
    
      +\displaystyle\int_{\mathbb{D}} |x|^k\frac{|U_{\varepsilon, s}(y)| |U_{\varepsilon, s}(y)-u_{\varepsilon,s}(y)|}{|x-y|^{n+2s}}dx dy. &
\end{array}
\end{equation}

We increase the two terms of the sum separately each one
\begin{equation}\label{20}
\begin{array}{ll}
 \displaystyle\int_{\mathbb{D}} |x|^k\frac{|U_{\varepsilon, s}(y)| |U_{\varepsilon, s}(y)-u_{\varepsilon,s}(y)|}{|x-y|^{n+2s}}dx dy &
 \leq 2 \displaystyle\int_{\mathbb{D}} |x|^k\frac{|U_{\varepsilon, s}(y)|^2}{|x-y|^{n+2s}}dx dy \\
 \\
 & \leq C \varepsilon^{n-2s} \displaystyle\int_{x\in B(0, \eta), y\in B(0, \eta)^c,  |x-y|>\frac{\eta}{2}}\frac{|x|^k}{|x-y|^{n+2s}} \\
 \\
  &\leq C \varepsilon^{n-2s} \eta^k \displaystyle\int_{x\in B(0, \eta), y\in B(0, \eta)^c,  |x-y|>\frac{\eta}{2}}\frac{1}{|x-y|^{n+2s}} dx dy\\
  \\
 & =\mathcal{O}(\varepsilon^{n-2s}),

\end{array}
\end{equation}

as $\varepsilon$ tends to zero,  and we have by recalling \eqref{fonction test u epsilon} and \eqref{majoration fonction test u} and as in (14.30) of \cite{BRS} 

\begin{equation}\label{21}
\begin{array}{ll}
 \displaystyle\int_{\mathbb{D}} |x|^k\frac{|U_{\varepsilon, s}(x)| |U_{\varepsilon, s}(y)-u_{\varepsilon,s}(y)|}{|x-y|^{n+2s}}dx dy &
 \leq  \displaystyle\int_{\mathbb{D}} |x|^k\frac{|U_{\varepsilon, s}(x)| (|U_{\varepsilon, s}(y)|)+(|u_{\varepsilon,s}(y)|)}{|x-y|^{n+2s}}dx dy \\
 \\
  &\leq 2 \displaystyle\int_{\mathbb{D}} |x|^k\frac{|U_{\varepsilon, s}(x)| |U_{\varepsilon, s}(y)|}{|x-y|^{n+2s}}dx dy \\
 \\
  &= \mathcal{O}(\varepsilon^{n-2s}),

\end{array}
\end{equation}
thus, by \eqref{20} and \eqref{21} together with \eqref{23} we get
\begin{equation}\label{25}
    \displaystyle\int_{\mathbb{D}} |x|^k\frac{|U_{\varepsilon, s}(x)-U_{\varepsilon, s}(y)| |U_{\varepsilon, s}(y)-u_{\varepsilon,s}(y)|}{|x-y|^{n+2s}}dx dy  \leq \mathcal{O}(\varepsilon^{n-2s})
\end{equation}
Therefore, using \eqref{2ème D} and \eqref{25} together with \eqref{somme D} we get
\begin{equation}\label{30}
      \displaystyle\int_{\mathbb{D}} |x|^k\frac{|u_{\varepsilon,s}(x)-u_{\varepsilon,s}(y)|^2}{|x-y|^{n+2s}}dx dy \leq  \displaystyle\int_{\mathbb{D}} |x|^k\frac{|(U_{\varepsilon, s}(x)-U_{\varepsilon, s}(y))|^2}{|x-y|^{n+2s}}dx dy + \mathcal{O}(\varepsilon^{n-2s})
\end{equation}

So by  Proposition \ref{As,k}, \eqref{45}, \eqref{40} and \eqref{30} we have
\begin{equation}\begin{array}{ll}
\displaystyle\int_{\mathbb{R}^n}\displaystyle\int_ {\mathbb{R}^n} (p(x)-p_0)\displaystyle\frac{|u_{\varepsilon,s}(x)-u_{\varepsilon,s}(y)|^2}{|x-y|^{n+2s}}dx dy & \leq \kappa C\varepsilon^{2s} + 2 \kappa \displaystyle\int_{\mathbb{D}} |x|^k\frac{|U_{\varepsilon, s}(x)-U_{\varepsilon, s}(y)|^2}{|x-y|^{n+2s}}dx dy \\
\\
& + \mathcal{O}(\varepsilon^{n-2s}). 
\end{array}
\end{equation}
Therefore, by this and \eqref{10}, for $n$ and $s$ as in \eqref{conditions n s} (here we have already taken into account the conditions in \eqref{contrainte 1}) and for $2\leq k< n-4s$ we write
\begin{equation}\label{première exp}
\begin{array}{ll}
\displaystyle\int_{\mathbb{R}^n}\displaystyle\int_ {\mathbb{R}^n} p(x)\frac{|u_{\varepsilon,s}(x)-u_{\varepsilon,s}(y)|^2}{|x-y|^{n+2s}}dx dy& \leq p_0 K_s+\kappa C\varepsilon^{2s}+ \mathcal{O}(\varepsilon^{n-2s})
\\
\\
&+ 2 \kappa \displaystyle\int_{\mathbb{D}} |x|^k\frac{|U_{\varepsilon, s}(x)-U_{\varepsilon, s}(y)|^2}{|x-y|^{n+2s}}dx dy. 
\end{array}
\end{equation}

Now,  we prove that the last term in the sum on the right hand-side of the previous inequality is $\mathcal{O}(\varepsilon^{k+2s})+\mathcal{O}(\varepsilon^{n-2s}).$\\
Note that for $(x,y) \in \mathbb{D}$ we have $x \in B(0,\eta)$ and $y \in B(0,\eta)^c$ 
\begin{equation}\label{12}
\begin{array}{ll}
\displaystyle\int_{\mathbb{D}} |x|^k\frac{|U_{\varepsilon, s}(x)-U_{\varepsilon, s}(y)|^2}{|x-y|^{n+2s}}dx dy \\
\\
\leq  \displaystyle\int_{B(0, \eta)} |x|^k \displaystyle\int_{B(0, \eta)^c}\frac{|U_{\varepsilon, s}(x)|^2+|U_{\varepsilon, s}(y)|^2}{|x-y|^{n+2s}}1_{B(0,\frac{\eta}{2})^c} (x-y)dx dy \\
\\
 \leq  \displaystyle\int_{B(0, \eta)} |x|^k \displaystyle\int_{B(0, \eta)^c}\frac{|U_{\varepsilon, s}(x)|^2}{|x-y|^{n+2s}}1_{B(0,\frac{\eta}{2})^c}(x-y)  dx dy
 \\
\\ 
 +  C \displaystyle\int_{B(0, \eta)} |x|^k \displaystyle\int_{B(0, \eta)^c} |U_{\varepsilon, s}(y)|^2 1_{B(0,\frac{\eta}{2})^c}(x-y)  dx dy. 
\end{array}
\end{equation}

For the second term on the right hand-side of the previous inequality, by \eqref{fonction test u} and the dominated convergence theorem as $\varepsilon$ tends to zero, and by \eqref{conditions n s} we have
\begin{equation}\label{conditions 0}
   \hbox{ for }  n\geq 4    \hbox{ and } 0<s<1
    \hbox{ or } n=3 \hbox{ and } 0<s<\frac{3}{4}
\end{equation}

then
\begin{equation}\label{26}
\begin{array}{ll}
   \displaystyle\int_{B(0, \eta)} |x|^k \displaystyle\int_{B(0, \eta)^c} |U_{\varepsilon, s}(y)|^2 1_{B(0,\frac{\eta}{2})^c}(x-y)  dx dy &
\leq \displaystyle\int_{B(0, \eta)} |x|^k dx\displaystyle\int_{B(0, \eta)^c} |U_{\varepsilon, s}(y)|^2   dy\\
\\
& \leq C_1 \displaystyle\int_{B(0, \eta)^c} |U_{\varepsilon, s}(y)|^2   dy \\
\\
& = C_1 \varepsilon^{n-2s} \bigg[ \displaystyle\int_{|y|\geq \eta} \displaystyle\frac{dy}{|y|^{2(n-2s)}}+o(1) \bigg]
\\
\\
  &  = \mathcal{O}(\varepsilon^{n-2s}),
\end{array}
\end{equation}
where $C$ and $C_1$ are positive constants possibly depending on $n,s,k$ and $\eta.$
So,  it remains to prove that 
\begin{equation}\label{11}
 \displaystyle\int_{B(0, \eta)} |x|^k \displaystyle\int_{B(0, \eta)^c}\frac{|U_{\varepsilon, s}(x)|^2}{|x-y|^{n+2s}}1_{B(0,\frac{\eta}{2})^c}(x-y)dy dx = \mathcal{O}(\varepsilon^{k+2s}).
 \end{equation}
In fact,  let $g_{1, k, \varepsilon}$ and $h$ be two functions such that,  for all $x\in\mathbb{R}^n$,  $g_{1, k, \varepsilon}(x) :=|x|^k |U_{\varepsilon, s}(x)|^2$ and for all $z\in\mathbb{R}^n$,  $h(z) :=\frac{1_{B(0,\frac{\eta}{2})^c}(z)}{|z|^{n+2s}}$. Since $s>0$, $h\in L^1(\mathbb{R}^n)$ .\\
We have  
\begin{equation}
\begin{array}{ll}
 ||g_{1, k, \varepsilon}||_{L^1(\mathbb{R}^n)} & =\displaystyle\int_{\mathbb{R}^n} |x|^k |U_{\varepsilon, s}(x)|^2 dx  \\
 \\
 & =  \varepsilon^{n-2s} \displaystyle\int_{\mathbb{R}^n} \frac{|x|^k}{(\varepsilon^2+|x|^2)^{n-2s}} dx \\
 \\
& = \varepsilon^{-n+2s} \displaystyle\int_{\mathbb{R}^n}\displaystyle\frac{|x|^k}{ \big( 1+\frac{|x|^2}{\varepsilon^2}\big)^{n-2s}}dx \\
\\
& =  \varepsilon^{k+2s} \displaystyle\int_{\mathbb{R}^n}\frac{|z|^k}{(1+|z|^2)^{n-2s}}dz. 
\end{array}
\end{equation}
Since $k<n-4s$ (for $n$ and $s$ as in \eqref{conditions n s}),  the integral $\displaystyle\int_{\mathbb{R}^n}\frac{|z|^k}{(1+|z|^2)^{n-2s}}dz$ is finite. Therefore,  $||g_{1,k,\varepsilon}||_{L^1(\mathbb{R}^n)} = \mathcal{O}(\varepsilon^{k+2s})$ and we have
\begin{equation}
\begin{array}{ll}
 \displaystyle\int_{B(0, \eta)} |x|^k \displaystyle\int_{B(0, \eta)^c}\frac{|U_{\varepsilon, s}(x)|^2}{|x-y|^{n+2s}}1_{B(0,\frac{\eta}{2})^c}(x-y)dx dy  
 \\
 \\
 \leq \displaystyle\int_{\mathbb{R}^n} \bigg( \displaystyle\int_{\mathbb{R}^n}|x|^k |U_{\varepsilon, s}(x)|^2\frac{1_{B(0,\frac{\eta}{2})^c}(x-y)}{|x-y|^{n+2s}}dx\bigg) dy
 \\
 \\
 = ||g_{1,k,\varepsilon}*h||_{L^1(\mathbb{R}^n)}  \\
 \\
  \leq ||g_{1,k,\varepsilon}||_{L^1(\mathbb{R}^n)} ||h||_{L^1(\mathbb{R}^n)} \\
 \\
  = \mathcal{O}(\varepsilon^{k+2s}),
\end{array}
\end{equation}
then, \eqref{11} yields.\\
Thus, by \eqref{26} and \eqref{11} together with \eqref{12}, for $n$ and $s$ as described in \eqref{conditions n s} (where conditions on $n$ and $s$ \eqref{conditions 0} are satisfied) and for $2\leq k<n-4s,$ 
\begin{equation}\label{dernière exp}
    \begin{array}{ll}
     \displaystyle\int_{\mathbb{D}} |x|^k\frac{|U_{\varepsilon, s}(x)-U_{\varepsilon, s}(y)|^2}{|x-y|^{n+2s}}dx dy & 
     \leq \mathcal{O}(\varepsilon^{n-2s}) + \mathcal{O}(\varepsilon^{k+2s}).
    \end{array}
\end{equation}

Finally, thanks to \eqref{première exp} and \eqref{dernière exp}, for $n$ and $s$ as described in \eqref{conditions n s} and for $2\leq k<n-4s$
we get the assertion \eqref{inégalité clé}
\begin{equation*}
 \displaystyle\int_{\mathbb{R}^n}\displaystyle\int_ {\mathbb{R}^n} p(x)  \frac{|u_{\varepsilon, s, a}(x)-u_{\varepsilon, s, a}(y)|^2}{|x-y|^{n+2s}}dx dy\leq p_0 K_s +\kappa C\varepsilon^{2s}+ \mathcal{O}(\varepsilon^{n-2s}) + \mathcal{O}(\varepsilon^{k+2s}), 
 \end{equation*}
which completes the proof of Theorem \ref{inégalité p S}.
\end{preuve}

\begin{remarque}
    Note that if we take $v_{\varepsilon,s,a}:={u_{\varepsilon,s,a}\over \vert\vert u_{\varepsilon,s,a}\vert\vert_{L^{q_s}(\mathbb{R}^n)}} $ inequality \eqref{inégalité clé} becomes 
    \begin{equation}
        \displaystyle\int_{\mathbb{R}^n}\displaystyle\int_ {\mathbb{R}^n} p(x)  \frac{|v_{\varepsilon, s, a}(x)-v_{\varepsilon, s, a}(y)|^2}{|x-y|^{n+2s}}dx dy\leq p_0 S_s +\kappa C\varepsilon^{2s}+ \mathcal{O}(\varepsilon^{n-2s}) + \mathcal{O}(\varepsilon^{k+2s}).
    \end{equation}
\end{remarque}

We conclude the proof of  Theorem \ref{théorème d'existence}. 
Thanks to Theorem \ref{inégalité p S}, \eqref{norme qs de u epsilon} and \eqref{norme 2 u epsilon} we write the inequality satisfied by the energy
\begin{equation}\label{Elambda}
\begin{array}{ll}
 S_{s,\lambda}(p) & \leq  E_\lambda(v_{\varepsilon,s,a}) \\
\\
& = E_\lambda \bigg(\displaystyle\frac{u_{\varepsilon,s,a}}{ \vert\vert u_{\varepsilon,s,a}\vert\vert_{L^{q_s}(\mathbb{R}^n)}}\bigg) 
 \\
 \\
 & \leq  p_0 S_s +\kappa C\varepsilon^{2s}- \lambda \frac{K_{2,s}}{K_{q_s}^{\frac{2}{q_s}}} \varepsilon^{2s}+ \mathcal{O}(\varepsilon^{n-2s})+\mathcal{O}(\varepsilon^{k+2s}) \\
 \\
 & = p_0 S_s + \varepsilon^{2s}\bigg(-\lambda \frac{K_{2,s}}{K_{q_s}^{\frac{2}{q_s}}}+\kappa C +\mathcal{O}(\varepsilon^{n-4s})+\mathcal{O}(\varepsilon^k)\bigg). \\
 \\
\end{array}
\end{equation}
We deduce that if  $n =  3$ and $s\in ]0, \frac{1}{4}]$ or $n =  4$ and $ s\in ]0,\frac{1}{2}]$ or  $n =  5$ and $s\in ]0, \frac{3}{4}]$ or 
if $n\geq 6$ and $s\in ]0,1[$ and for $k\in[2,n-4s[,$ there exists a constant $C =C( n,s,k) > 0$ such that  for every $\kappa \in ]0, C\lambda[$ we have  $S_{s, \lambda}(p)\leq E_\lambda(v_{\varepsilon,s,a})<p_0 S_s$. 
In all the cases above $S_{s, \lambda}(p)$ is achieved for $0<\lambda<\lambda_{1, p, s}$ and for every $\kappa \in ]0,C\lambda[$  thus the Theorem \ref{théorème d'existence} follows at once. \\
\\
We end this section by stating the following remark which  ensure the existence of non negative solution for our minimization problem.
\begin{remarque} \label{solution positive1}  
\begin{enumerate}
    \item As a consequence of  Remark \ref{solution positive} and the proof of Theorem \ref{théorème d'existence} we obtain actually that $S_{s, \lambda}(p)$ possesses a positive solution.
    \item Let $u$ be a positive function which achieves the infimum $S_{s, \lambda}(p)$. Therefore, for all $\phi \in \mathbb{H}_0^s(\Omega)$, the function $u$ satisfies the following Euler's equation
\begin{equation}\label{theta}
\begin{array}{ll}
        \displaystyle\int_{\mathbb{R}^n\times\mathbb{R}^n} p(x)\displaystyle\frac{\big(u(x)-u(y)\big) \big(\phi(x)-\phi(y)\big)}{|x-y|^{n+2s}}dy dx -  \lambda \displaystyle\int_\Omega u(x) \phi(x) dx \\
        \\= S_{s,\lambda}(p) \displaystyle\int_\Omega |u(x)|^{q_s-2} u(x) \phi(x) dx.
        \end{array}
\end{equation}
   Note that since $\lambda<\lambda_{1,p,s}$,  $S_{s,\lambda}(p)>0$. It follows that, $v:=(S_{s,\lambda}(p))^\frac{1}{2-q_s}u$ is positive a solution for the following variational problem 
    \begin{equation}\label{prob sans theta}
        \displaystyle\int_{\mathbb{R}^n\times\mathbb{R}^n} p(x)\displaystyle\frac{\big(v(x)-v(y)\big) \big(\phi(x)-\phi(y)\big)}{|x-y|^{n+2s}}dy dx -  \lambda \displaystyle\int_\Omega u(x) \phi(x) dx = \displaystyle\int_\Omega |v(x)|^{q_s-2} v(x) \phi(x) dx,
\end{equation}
 $\phi \in \mathbb{H}_0^s(\Omega)$.
Let $\Phi_{1, p, s}$ an eigenvector associated to $\lambda_{1, p, s}$. We take $\phi=\phi_{1, p, s}$,  we get
\begin{equation}
\begin{array}{ll}
\displaystyle\int_{\mathbb{R}^n}\displaystyle\int_{\mathbb{R}^n} p(x) \frac{(v(x)-v(y) (\phi_{1, p, s}(x)-\phi_{1, p, s}(y))}{|x-y|^{n+2s}}dxdy - \lambda \displaystyle\int_\Omega v(x) \phi_{1, p, s}(x) dx \\
\\= \displaystyle\int_\Omega |v(x)|^{q_s-2} v(x) \phi_{1, p, s}(x) dx, 
\end{array}
\end{equation}
and
\begin{equation}
\begin{array}{ll}
(\lambda_{1, p, s}-\lambda) \displaystyle\int_{\Omega}   v(x) \Phi_{1, p, s}(x) dx = \displaystyle\int_{\Omega}|v(x)|^{q_s-2} v(x) \phi_{1, p, s}(x) dx.
\end{array}
\end{equation}
Identically as in (\cite{BRS}, chapter3),   we can prove that  $\phi_{1, p, s}>0$.  Therefore,  we deduce 
\begin{equation}
\displaystyle\int_{\Omega}|v|^{q_s-2}\phi_{1, p, s} dx>0.
\end{equation}
Then we get
$(\tilde\lambda_{1, p, s}-\lambda)\displaystyle\int_{\mathbb{R}^n} v(x) \phi_{1, p, s}(x) dx > 0$ 
and so $\tilde\lambda_{1, p, s}>\lambda$.\\
Finally, if $\lambda \geq \tilde\lambda_{1, p, s} $, there exists no solution achieving the infimum \eqref{infimum} and problem \eqref{prob sans theta} does not have any solution.  
\end{enumerate}
\end{remarque}
\section{Existence of non-ground state solutions}

Let $p$ be as described in \eqref{def de p} such that $p$ is bounded. Let $n\geq3,$ $s\in]0,1[,$ $\lambda>0$ and $q\in[2,q_s[,$ with $q_s=\frac{2n}{n-2s}$.

It's obvious to see that equation \eqref{formulation faible} is the Euler-Lagrange equation of the functional $\Phi_{p,s}: \mathbb{H}_0^s(\Omega) \rightarrow \R$ defined as 

\begin{equation}\label{functional}
    \Phi_{p,s}(u) :=\frac{1}{2} \displaystyle\int_{\mathbb{R}^n \times \mathbb{R}^n} p(x) \frac{|u(x)-u(y)|^2}{|x-y|^{n+2s}}dx dy - \frac{\lambda}{q} \displaystyle\int_{\Omega}|u(x)|^q dx -\frac{1}{q_s} \displaystyle\int_{\Omega}|u(x)|^{q_s}dx.
\end{equation}
Notice that the functional $\Phi_{p,s}$ does not satisfy the Palais-Smale condition globally and this is due to the lack of compactness of the embedding $\mathbb{H}_0^s(\Omega)$ into $L^{q_s}(\Omega)$. As a consequence, an estimate of the critical level of $\Phi_{p,s}$ is necessary and as we will see in what follows, we cannot apply the classical mountain pass Theorem, that's why we will apply a variant of the mountain pass Theorem without the Palais-Smale condition, as given in   (\cite{BN1},Theorem 2.2).\\
Thanks to the fact that the embedding $\mathbb{H}^s_0(\Omega)$ into $L^r(\Omega)$ is compact for $r \in [1,q_s[$ and it is continuous for $r=q_s$ with $\Omega$ a domain with continuous boundary, the energy functional $\Phi_{p,s}$ is well defined.\\
Moreover, $\Phi_{p,s}$ is Fréchet differentiable in $u\in \mathbb{H}^s_0(\Omega)$, and, for any $\varphi \in\mathbb{H}^s_0(\Omega) $ and $ q \in [2,q_s[$,
\begin{equation}
    \begin{array}{ll}
    <\Phi'_{p,s}(u),\varphi> & = \displaystyle\int_{\mathbb{R}^n\times \mathbb{R}^n} p(x)\frac{(u(x)-u(y))(\varphi(x)-\varphi(y))}{|x-y|^{n+2s}}dxdy \\
    \\
    & -\lambda \displaystyle\int_{\Omega} |u(x)|^{q-2} u(x) \varphi(x) dx  - \displaystyle\int_{\Omega} |u(x)|^{q_s-2}u(x)\varphi(x) dx .
    \end{array}
\end{equation}
Notice that critical points of $\Phi_{p,s}$ are solutions to problem \eqref{formulation faible}. In order to find these critical points, we proceed as follows; first we start proving that  $\Phi_{p,s}$ has a suitable geometric sturcture as stated in conditions (2.9) and (2.10) of (\cite{BN1},Theorem 2.2). 

\begin{prop}\label{vérif 2.9}
 1) Let $q \in ]2,q_s[$, $\lambda>0$. There exist $\rho >0$ and $\beta >0$ such that, for all $u\in\mathbb{H}^s_0(\Omega)$ with $\|u\|_ {\mathbb{H}^s_0(\Omega)}=\rho$, we have $\Phi_{p,s}(u)\geq \beta$.\\
 \\
 2) For $q=2$ and $\lambda\in]0,\lambda_{1,p,s}[$, There exist $\rho >0$ and $\beta >0$ such that, for all $u\in\mathbb{H}^s_0(\Omega)$ with $\|u\|_ {\mathbb{H}^s_0(\Omega)}=\rho$, we have $\Phi_{p,s}(u)\geq \beta$.
\end{prop}

\begin{preuve}
    Let  $ q \in ]2,q_s[$ and $u\in \mathbb{H}^s_0(\Omega)$ such that $ \|u\|_{ \mathbb{H}^s_0(\Omega)}\leq 1$. Thanks to the continuous Sobolev embeddings $\mathbb{H}^s_0(\Omega)$ into $L^{r}(\Omega),$ where $1\leq r\leq q_s$, we have
    \begin{equation}\label{minoration 2.9}
        \begin{array}{ll}
          \Phi_{p,s}(u) & =\frac{1}{2} \displaystyle\int_{\mathbb{R}^n \times \mathbb{R}^n} p(x) \frac{|u(x)-u(y)|^2}{|x-y|^{n+2s}}dx dy - \frac{\lambda}{q} \displaystyle\int_{\Omega}|u(x)|^q dx -\frac{1}{q_s} \displaystyle\int_{\Omega}|u(x)|^{q_s}dx \\
          \\
          & \geq\frac{1}{2}p_0 \|u\|_{ \mathbb{H}^s_0(\Omega)}^2 -\frac{\lambda }{q} \alpha_q^{-\frac{q}{2}} \|u\|_{ \mathbb{H}^s_0(\Omega)}^q -\frac{1}{q_s} S_s^{-\frac{q_s}{2}} \|u\|_{ \mathbb{H}^s_0(\Omega)}^{q_s} \\
          \\
          &= \|u\|_{ \mathbb{H}^s_0(\Omega)}^2 \bigg( \frac{1}{2}p_0-\frac{\lambda}{q}\alpha_q^{-\frac{q}{2}}\|u\|_{ \mathbb{H}^s_0(\Omega)}^{q-2}-\frac{1}{q_s} S_s^{-\frac{q_s}{2}} \|u\|_{ \mathbb{H}^s_0(\Omega)}^{q_s-2}\bigg),
        \end{array}
    \end{equation}
    where $\alpha_q$ is the best Sobolev constant of the embedding $\mathbb{H}^s_0(\Omega)$ into $L^{q}(\Omega).$\\
Now let $u\in \|u\|_{ \mathbb{H}^s_0(\Omega)}$ be such that $ \|u\|_{ \mathbb{H}^s_0(\Omega)}= \rho >0$. Since $q_s > 2$, we choose $\rho$ sufficiently small in order to get $ \big( \frac{1}{2}p_0-\frac{\lambda}{q}\alpha_q^{-\frac{q}{2}}\rho^{q-2}-\frac{1}{q_s} S_s^{-\frac{q_s}{2}} \rho^{q_s-2}\big)>0$ and by defining $\beta:= \rho^2 \big( \frac{1}{2}p_0-\frac{\lambda}{q}\alpha_q^{-\frac{q}{2}}\rho^{q-2}-\frac{1}{q_s} S_s^{-\frac{q_s}{2}} \rho^{q_s-2}\big) >0 $, we obtain  $\Phi_{p,s}(u)\geq \beta >0$.

    For the case where $q= 2$, let $0<\lambda <\lambda_{1,p,s}$. We recall that the norm associated with the scalar product \eqref{produit scalaire p} is equivalent to the ordinary norm over $\mathbb{H}_0^s(\Omega)$,  then we have
    \begin{equation}\label{minoration 2.9}
        \begin{array}{ll}
          \Phi_{p,s}(u) & =\frac{1}{2} \displaystyle\int_{\mathbb{R}^n \times \mathbb{R}^n} p(x) \frac{|u(x)-u(y)|^2}{|x-y|^{n+2s}}dx dy - \frac{\lambda}{2} \displaystyle\int_{\Omega}|u(x)|^2 dx -\frac{1}{q_s} \displaystyle\int_{\Omega}|u(x)|^{q_s}dx \\
          \\
          & \geq\frac{1}{2}\bigg( \|u\|_{ \mathbb{H}^s_0(\Omega)}^2 -\frac{\lambda }{\lambda_{1,p,s}} \|u\|_{ \mathbb{H}^s_0(\Omega)}^2\bigg) -\frac{1}{q_s} S_s^{-\frac{q_s}{2}} \|u\|_{ \mathbb{H}^s_0(\Omega)}^{q_s} \\
          \\
          &  \geq  \|u\|_{ \mathbb{H}^s_0(\Omega)}^2 \bigg[\frac{1}{2}\bigg(1-\frac{\lambda}{\lambda_{1,p,s}} \bigg)- \frac{1}{q_s} S_s^{-\frac{q_s}{2}}\|u\|_{ \mathbb{H}^s_0(\Omega)}^{q_s-2}\bigg].
        \end{array}
    \end{equation}
Hence, it easily follows that for $\lambda \in ]0,\lambda_{1,p,s}[$,
\begin{equation}
    \Phi_{p,s}(u) \geq \tilde{\alpha}_1  \|u\|_{ \mathbb{H}^s_0(\Omega)}^2 \big(1- \tilde{\alpha}_2 \|u\|_{ \mathbb{H}^s_0(\Omega)}^{q_s-2}\big),
\end{equation}
for suitable positive constants $\tilde{\alpha}_1$ and $\tilde{\alpha}_2$.
Now let $u\in \|u\|_{ \mathbb{H}^s_0(\Omega)}$ be such that $ \|u\|_{ \mathbb{H}^s_0(\Omega)}= \rho >0$. Since $q_s > 2$, we choose $\rho$ sufficiently small in order to get $1-\tilde{\alpha}_2 \rho^{q_s-2} >0$ and by defining $\beta:= \tilde{\alpha}_1 \rho^2 (1-\tilde{\alpha}_2\rho^{q_s-2}) >0 $, we obtain  $\Phi_{p,s}(u)\geq \beta >0$.
\end{preuve}

\begin{prop}\label{vérif 2.10'}
   Let $ q \in [2,q_s[$ and $\lambda>0$. There exists $e\in \mathbb{H}^s_0(\Omega) $ such that $e\geq 0$ a.e in $\mathbb{R}^n$, $\|e\|_{\mathbb{H}^s_0(\Omega)}> \rho$, and $\Phi_{p,s}(e)<\beta,$ where $\rho$ and $\beta$ are given in Proposition \ref{vérif 2.9}.
\end{prop}
\begin{preuve}
    Let $u$ be a fixed function in $\mathbb{H}^s_0(\Omega)$, such that $\|u\|_{\mathbb{H}^s_0(\Omega)}=1$ and $u\geq 0$ a.e in $\mathbb{R}^n$; we remark that this choice is possible. In fact, we replace any $u\in \mathbb{H}^s_0(\Omega)$ by its positive part since, if $u\in \mathbb{H}^s_0(\Omega)$, $u^+ \in \mathbb{H}^s_0(\Omega)$ where $u^+(x):=max\{ u(x),0\}.$\\
Let $\zeta >0$, since $p$ is bounded, there exists a positive constant $M$ such that for all $x \in \mathbb{R}^n$, $p(x) \leq M$, and we have

\begin{equation}
\begin{array}{ll}
         \Phi_{p,s}(\zeta u) & =\displaystyle\frac{\zeta^2}{2}\displaystyle\int_{\mathbb{R}^n\times\mathbb{R}^n} p(x) \displaystyle
         \frac{|u(x)-u(y)|^2}{|x-y|^{n+2s}}dx dy  \\
         \\
         & -\displaystyle\frac{\lambda}{q} \zeta^q \displaystyle\int_{\Omega}|u(x)|^q dx -\displaystyle\frac{\zeta^{q_s}}{q_s} \displaystyle\int_{\Omega}|u(x)|^{q_s}dx\\
         \\
         & \leq\displaystyle\frac{\zeta^2}{2} M \| u \|_ {\mathbb{H}_0^s(\Omega)}^2-\displaystyle\frac{\lambda}{q} \zeta^q \displaystyle\int_{\Omega} |u(x)|^q dx -\displaystyle\frac{\zeta^{q_s}}{q_s} \displaystyle\int_\Omega |u(x)|^{q_s}dx \\
         \\
         & =\displaystyle\frac{M \zeta^2}{2}  -\displaystyle\frac{\lambda}{q} \zeta^q \displaystyle\int_\Omega |u(x)|^q dx -\displaystyle\frac{\zeta^{q_s}}{q_s}\displaystyle\int_\Omega |u(x)|^{q_s} dx
         \\
         \\
         & \leq \displaystyle\frac{M \zeta^2}{2} -\displaystyle\frac{\zeta^{q_s}}{q_s}\displaystyle\int_\Omega |u(x)|^{q_s} dx.
\end{array}
\end{equation}
     Since $q_s>2$, by passing to the limit $\zeta \rightarrow +\infty$ and for $u$ fixed, we get $$\Phi_{p,s}(\zeta u) \rightarrow -\infty,$$
 so we take $e=\zeta u$, with $\zeta$ large enough so that we get $\Phi_{p,s}(\zeta u)<0$ and the assertion yields since $\beta>0.$
\end{preuve}
\begin{prop}\label{vérif 2.11}
 Let $n\geq 3$, $q\in [2,q_s[$ and $s\in ]0,1[$. 
 Let 
    \begin{equation}
    c := \inf_{P\in \mathcal{P}}\sup_{v\in P([0,1])}\Phi_{p,s}(v),
    \end{equation}
    where
    \begin{equation}
    \mathcal{P} := \{ P\in C([0,1];\mathbb{H}^s_0(\Omega)): P(0)=0, P(1) = e    \},
    \end{equation}
    with $e$ given in Proposition \ref{vérif 2.10'}. Then $\beta \leq c<\frac{s}{n}(p_0 S_s)^\frac{n}{2s},$ where $\beta$ given 
    in Proposition \ref{vérif 2.9}. 
\end{prop}
\begin{preuve}
Let $\rho$ as given in Proposition \ref{vérif 2.9}.
    For all $P \in \mathcal{P},$ the function $t \mapsto \|P(t) \|_{\mathbb{H}_0^s(\Omega)}$ is continuous in $[0,1]$, $\|P(0)\|_{\mathbb{H}_0^s(\Omega)} <\rho$ and $\|P(1)\|_{\mathbb{H}_0^s(\Omega)}=\|e\|_{\mathbb{H}_0^s(\Omega)}>\rho$, therefore, there exists a real $\bar{t} \in ]0,1[$ such that $\| P(\bar{t})\|= \rho.$  Then, by Proposition \ref{vérif 2.9} and taking $P(\bar{t})$ as a function test 
    \begin{equation}
    \max_{v\in P([0,1])}\Phi_{p,s}(v)\geq \Phi_{p,s}(P(\bar{t}))\geq \displaystyle\inf_{\substack{v \in \mathbb{H}_0^s(\Omega) \\ \|v\|_{\mathbb{H}_0^s(\Omega)}=\rho}}\Phi_{p,s}(v) \geq \beta
      \end{equation}
      
Finally, we obtain $c\geq \beta$.\\
Now Let us prove that $c<\frac{s}{n}(p_0 S_s)^\frac{n}{2s}$ and let $v_{\varepsilon,s,a} := \displaystyle\frac{u_{\varepsilon,s,a}}{||u_{\varepsilon,s,a}||_{L^{q_s}(\mathbb{R}^n)}}$ where $u_{\varepsilon,s,a}$ is defined in \eqref{fonction test u epsilon}.\\
\\
Let us consider $\tilde{X}_\varepsilon := \displaystyle\int_{\mathbb{R}^n} \displaystyle\int_{\mathbb{R}^n} p(x) \frac{|v_{\varepsilon,s,a}(x)-v_{\varepsilon,s,a}(y)|^2}{|x-y|^{n+2s}} dy dx$, then we write for $t\geq 0$
\begin{equation*}
\Phi_{p,s}(tv_{\varepsilon,s,a})=\frac{1}{2} t^2 \tilde{X}_\varepsilon - \frac{t^{q_s}}{q_s}-\frac{t^q}{q} \lambda \displaystyle\int_\Omega |v_{\varepsilon,s,a}|^q dx \leq\frac{1}{2} t^2 \tilde{X}_\varepsilon - \frac{t^{q_s}}{q_s},
\end{equation*}
therefore $\displaystyle\lim_{t \rightarrow +\infty} \Phi_{p,s}(t v_{\varepsilon,s,a})=-\infty$ and $\displaystyle\sup_{t\geq 0} \Phi_{p,s}(t v_{\varepsilon,s,a})$ is achieved at some $t_\varepsilon \geq 0$. If $t_\varepsilon = 0,$ then $\displaystyle\sup_{t\geq 0} \Phi_{p,s}(t v_{\varepsilon,s,a})=0$ and there is nothing to prove. So we assume that $t_\varepsilon >0$. Since $\displaystyle\frac{\partial \Phi_{p,s}}{\partial t}(t_\varepsilon v_{\varepsilon,s,a}) = 0$, we have
\begin{equation}\label{dérivée s'annule}
    t_\varepsilon \tilde{X}_\varepsilon - t_\varepsilon^{q_s-1} -\lambda t_\varepsilon^{q-1}  \displaystyle\int_\Omega |v_{\varepsilon,s}(x)|^qdx=0,
\end{equation}
and so
\begin{equation}\label{majoration t epsilon}
t_\varepsilon \leq \tilde{X}_\varepsilon ^\frac{1}{q_s-2}.
\end{equation}
We set 
\begin{equation}
Y_\varepsilon := \displaystyle\sup_{t\geq0} \Phi_{p,s}(tv_{\varepsilon,s,a}) = \Phi_{p,s}(t_\varepsilon v_{\varepsilon,s,a}).
\end{equation}
The function $t \mapsto (\displaystyle\frac{1}{2} t^2 \tilde{X}_\varepsilon -\displaystyle\frac{t^{q_s}}{q_s})$ is increasing on the interval $[0, \tilde{X}_\varepsilon^\frac{1}{q_s-2}],$ we have, by \eqref{majoration t epsilon},
\begin{equation}
    \begin{array}{ll}
    Y_\varepsilon & =\displaystyle\frac{1}{2} t_\varepsilon^2 \tilde{X}_\varepsilon -\frac{t_\varepsilon^{q_s}}{q_s}-\displaystyle\frac{\lambda}{q} t_\varepsilon^{q}  \displaystyle\int_\Omega |v_{\varepsilon,s,a}(x)|^q dx \\
    \\
    & \leq\displaystyle\frac{1}{2} \tilde{X}_\varepsilon ^\frac{2}{q_s-2} \tilde{X}_\varepsilon -\frac{\tilde{X}_\varepsilon ^\frac{q_s}{q_s-2}}{q_s}-\displaystyle\frac{\lambda}{q} t_\varepsilon^{q}  \displaystyle\int_\Omega |v_{\varepsilon,s,a}(x)|^q dx \\
    \\
    & =\displaystyle\frac{1}{2} \tilde{X}_\varepsilon ^\frac{n}{2s} -\frac{\tilde{X}_\varepsilon ^\frac{n}{2s}}{q_s} -\displaystyle\frac{\lambda}{q} t_\varepsilon^{q}  \displaystyle\int_\Omega |v_{\varepsilon,s,a}(x)|^q dx \\
    \\
    & =\displaystyle\frac{s}{n} \tilde{X}_\varepsilon^\frac{n}{2s} -\displaystyle\frac{\lambda}{q} t_\varepsilon^{q}  \displaystyle\int_\Omega |v_{\varepsilon,s,a}(x)|^q dx
    \end{array}
\end{equation}
We claim that 
\begin{equation}\label{lim t epsilon}
 t_\varepsilon \rightarrow (p_0 S_s)^\frac{1}{q_s-2},
\end{equation}
as $\varepsilon \rightarrow 0$.\\
In fact, by \eqref{dérivée s'annule} we have 
\begin{equation}\label{t epsilon}
\tilde{X}_\varepsilon-t_\varepsilon^{q_s-2} -\lambda t_\varepsilon^{q-2}  \displaystyle\int_\Omega |v_{\varepsilon,s,a}(x)|^qdx=0.
\end{equation}

\begin{lemme}
   Let $q\in [2,q_s[$. Let $n\geq4$ and $s \in ]0,1[,$ or  $n=3$ and $s < \frac{3}{4}$. For any $n > 4 s$ we have
\begin{equation}\label{norme q u epsilon}
||u_{\varepsilon,s,a}||_{L^q(\mathbb{R}^n)}^q = K_{q,s} \varepsilon^{n-\frac{q(n-2s)}{2}}+\mathcal{O}(\varepsilon^{\frac{q(n-2s)}{2}})
\end{equation}
with $K_{2, s}$ is a positive  constant,
$K_{q,s}=\displaystyle\int_{\mathbb{R}^n}\frac{dy}{(1+|y|^2)^{\frac{q(n-2s)}{2}}}. $
\end{lemme}
\begin{preuve}
Let $0<s<1$ and $q\in [2,q_s[$ and $\Psi$ as in \eqref{fonction test u epsilon}
      \begin{equation*}
      \begin{array}{ll}
\displaystyle\int_{\mathbb{R}^n}|u_{\varepsilon,s,a}(x)|^qdx & =\displaystyle\int_{\mathbb{R}^n}\bigg( \frac{\varepsilon}{\varepsilon^2+|x-a|^2}\bigg)^{\frac{q(n-2s)}{2}} \Psi^q(x)dx \\
& = \varepsilon^{\frac{q(n-2s)}{2}}\displaystyle\int_{\mathbb{R}^n} \frac{\Psi^q(x)-1}{\big(\varepsilon^2+|x-a|^2\big)^{\frac{q(n-2s)}{2}}}dx+\int_{\mathbb{R}^n} \bigg( \frac{\varepsilon}{\varepsilon^2+|x-a|^2}\bigg) ^{\frac{q(n-2s)}{2}}dx
 \end{array}
   \end{equation*}
   We have
 $\displaystyle\int_{\mathbb{R}^n} \frac{\Psi^q(x)-1}{\big(\varepsilon^2+|x-a|^2\big)^{\frac{q(n-2s)}{2}}}dx = \displaystyle\int_{B(a,\eta)^c} \frac{\Psi^q(x)-1}{\big(\varepsilon^2+|x-a|^2\big)^{\frac{q(n-2s)}{2}}}dx $ is bounded since $n >4s$, $q\geq 2$ and
 ~\\
 
 $$\bigg|\displaystyle\int_{\mathbb{R}^n\setminus{\mathbf{B}(a,\eta)}} \frac{\Psi^q(x)-1}{\big(\varepsilon^2+|x-a|^2\big)^{\frac{q(n-2s)}{2}}}dx\bigg|\leq \displaystyle\int_{\mathbb{R}^n\setminus{\mathbf{B}(a,\eta)}} \frac{\|\Psi\|_{\mathbf{L}^\infty}^q+1}{|x-a|^{q(n-2s)}}dx<\infty.$$

We conclude then for $n>4s$, $q\geq 2$ 
$$\displaystyle\int_{\mathbb{R}^n}|u_{\varepsilon,s,a}(x)|^qdx=\varepsilon^{\frac{q(n-2s)}{2}} \mathcal{O}(1)+\displaystyle\int_{\mathbb{R}^n}\bigg( \frac{\varepsilon}{\varepsilon^2+|x-a|^2}\bigg) ^{\frac{q(n-2s)}{2}}dx.$$

Now, let us study the integral $\displaystyle\int_{\mathbb{R}^n}\bigg( \frac{\varepsilon}{\varepsilon^2+|x-a|^2}\bigg) ^{\frac{q(n-2s)}{2}}dx$. Proceeding to the changement of variables $\varepsilon y = x-a$, we obtain
\begin{equation*}   
\displaystyle\int_{\mathbb{R}^n}\bigg( \frac{\varepsilon}{\varepsilon^2+|x-a|^2}\bigg) ^{\frac{q(n-2s)}{2}}dx    
         =\varepsilon ^{n-{\frac{q(n-2s)}{2}}} \displaystyle\int_{\mathbb{R}^n} \frac{dy}{(1+|y|^2)^{\frac{q(n-2s)}{2}}},
\end{equation*}
 
if $n>4s$ and $q\geq 2$ then $\displaystyle\int_{\mathbb{R}^n} \frac{dy}{(1+|y|^2)^{\frac{q(n-2s)}{2}}}$ converges,  thus
\\
$$\displaystyle\int_{\mathbb{R}^n}\bigg( \frac{\varepsilon}{\varepsilon^2+|x-a|^2}\bigg) ^{\frac{q(n-2s)}{2}}dx=\varepsilon^{n-{\frac{q(n-2s)}{2}}} K_{q,s},$$
where $K_{q,s}=\displaystyle\int_{\mathbb{R}^n} \frac{dy}{(1+|y|^2)^{\frac{q(n-2s)}{2}}}.$
Therefore, equation \eqref{norme q u epsilon} follows.
\end{preuve}\\
\\
Since $q<q_s$ using \eqref{norme qs de u epsilon}, \eqref{norme q u epsilon} and  \eqref{t epsilon}  we get 
 \begin{equation*}
\tilde{X}_\varepsilon-t_\varepsilon^{q_s-2} + o( t_\varepsilon^{q-2} ) =0.
\end{equation*}
which implies the claim \eqref{lim t epsilon}.\\
Using Theorem \ref{inégalité p S} together with \eqref{norme q u epsilon} and \eqref{lim t epsilon}, for $k\in[2,n-4s[$ and $q\in ]2,q_s[$ we have
\begin{equation}\label{inégalité Y epsilon}
\begin{array}{ll}
Y_\varepsilon & \leq \displaystyle\frac{s}{n} \tilde{X}_\varepsilon^\frac{n}{2s} -\displaystyle\frac{\lambda}{q} t_\varepsilon^{q}  \displaystyle\int_\Omega |v_{\varepsilon,s,a}(x)|^q dx \\
\\
& \leq \displaystyle\frac{s}{n} (p_0 S_s)^\frac{n}{2s}+ \kappa C \varepsilon^{2s} + \mathcal{O}(\varepsilon^{n-2s}) + \mathcal{O}(\varepsilon^{k+2s})
-  C_2 \varepsilon^{n-\frac{q(n-2s)}{2}}+\mathcal{O}(\varepsilon^{\frac{q(n-2s)}{2}})   \\
\\
& \leq \displaystyle\frac{s}{n}  (p_0 S_s)^\frac{n}{2s} - C_2 \varepsilon^{n-\frac{q(n-2s)}{2}} \bigg(1 + C_4 \varepsilon^{2s-n+\frac{q(n-2s)}{2}} + \mathcal{O}\big(\varepsilon^{-2s+\frac{q(n-2s)}{2}}\big)\\
\\
& +\mathcal{O}(\varepsilon^{-n+\frac{q(n-2s)}{2}+k+2s}) +\mathcal{O}(\varepsilon^{-n+q(n-2s)}) \bigg)
\end{array}
\end{equation}
 
where $C_2>0 $.\\
Therefore, since $n > 4s$ and $q\in ]2,q_s[$ we obtain for $\varepsilon>0$ sufficiently small 
\begin{equation}
Y_\varepsilon <\frac{s}{n} p_0^\frac{n}{2s} S_s^\frac{n}{2s},
\end{equation}
For $q=2$, we write
\begin{equation}\label{inégalité Y epsilon pour q=2}
\begin{array}{ll}
Y_\varepsilon & \leq \displaystyle\frac{s}{n} \tilde{X}_\varepsilon^\frac{n}{2s} -\displaystyle\frac{\lambda}{q} t_\varepsilon^{q}  \displaystyle\int_\Omega |v_{\varepsilon,s,a}(x)|^2 dx \\
\\
& \leq \displaystyle\frac{s}{n} (p_0 S_s)^\frac{n}{2s}+ \kappa C \varepsilon^{2s} + \mathcal{O}(\varepsilon^{n-2s}) + \mathcal{O}(\varepsilon^{k+2s}) 
\\
\\
& - \lambda C_1 \varepsilon^{2s} + \mathcal{O}(\varepsilon^{n-2s})  
\\
\\
& \leq \displaystyle\frac{s}{n}  (p_0 S_s)^\frac{n}{2s} + \varepsilon^{2s} \bigg(\kappa C - \lambda C_1 + \mathcal{O}\big(\varepsilon^{n-4s}\big)\\
\\
& +\mathcal{O}(\varepsilon^k) +\mathcal{O}(\varepsilon^{n-4s}) \bigg)
\end{array}
\end{equation}
Thus, for $n > 4s$ and $q = 2$,$\lambda \in]0,\lambda_{1,p,s}[,$ $\kappa \in ]0,C\lambda[$ and for $\varepsilon>0$ sufficiently small we obtain
\begin{equation}
Y_\varepsilon <\frac{s}{n} p_0^\frac{n}{2s} S_s^\frac{n}{2s},
\end{equation}

When $n = 3$, using the computation in \eqref{inégalité Y epsilon} and \eqref{norme q u epsilon} it results that for $Y_\varepsilon <\frac{s}{3} p_0^\frac{3}{2s} S_s^\frac{3}{2s}$ for $s\in ]0,\frac{3}{4}[$.

hence the proof is complete.
\end{preuve}

\subsection{Proof of Theorem \ref{prob sans minimisation}}

Proposition \ref{vérif 2.9} and \ref{vérif 2.10'} give that the geometry of the variant of the mountain pass Theorem given in (\cite{BN1}, Theorem 2.2) is fulfilled by $\Phi_{p,s}.$ Furthermore, we easily have $\Phi_{p,s}(0)=0 < \beta,$ with $\beta$ given in Proposition  \ref{vérif 2.9} and by Proposition \ref{vérif 2.11}, we deduce that all hypotheses of (\cite{BN1}, Theorem 2.2) are satisfied, then, there is a sequence $(u_j)_{j\in\N}$ in $\mathbb{H}_0^s(\Omega)$ such that

\begin{equation}\label{conv c}
\Phi_{p,s}(u_j)\rightarrow c
\end{equation}
and
\begin{equation}\label{conv 0}
sup \bigg \{ |<\Phi'_{p,s}(u_j),\varphi>|: \varphi \in \mathbb{H}_0^s(\Omega), \| \varphi \|_ {\mathbb{H}_0^s(\Omega)}=1      \bigg \} \rightarrow 0
\end{equation}
as $j$ tends to $+\infty.$\\

In order to finish the proof of Theorem \ref{prob sans minimisation} we proceed as in \cite{BN1} and (\cite{BRS}, chapter 14) and we need to establish the following result
    \begin{lemme}
    The sequence $(u_j)_{j\in\N}$ is bounded in $\mathbb{H}_0^s(\Omega)$.
    \end{lemme}
\begin{preuve}
    For any $j\in \N$, by \eqref{conv c} and \eqref{conv 0}, it easily follows that there exists $k>0$ such that
\begin{equation}\label{phi k}
    |\Phi_{p,s}(u_j)| \leq k
\end{equation}
and
\begin{equation}\label{phi' k}
    |<\Phi'_{p,s}(u_j),\displaystyle\frac{u_j}{\|u_j\|_{\mathbb{H}_0^s(\Omega)}}>| \leq k
\end{equation}
By \eqref{phi k} and \eqref{phi' k}, we have
\begin{equation}\label{phi}
    \Phi_{p,s}(u_j)-\frac{1}{2}<\Phi'_{p,s}(u_j),u_j>  \leq k (1+\|u_j\|_{\mathbb{H}_0^s(\Omega)})
\end{equation}
Moreover, we have 
\begin{equation}
 \begin{array}{ll}
 \Phi_{p,s}(u_j)-\frac{1}{2}<\Phi'_{p,s}(u_j),u_j> & = (-\frac{1}{q_s}+\frac{1}{2})\|u_j \|_{L^{q_s}(\Omega)}^{q_s}+\lambda(-\frac{1}{q}+\frac{1}{2}) \| u_j\|_{L^q(\Omega)}^q  \\
 \\
 & =\frac{s}{n}\|u_j \|_{L^{q_s}(\Omega)}^{q_s}+\lambda(-\frac{1}{q}+\frac{1}{2}) \| u_j\|_{L^q(\Omega)}^q,
 \end{array}
\end{equation}
then, by this and \eqref{phi}, we get that, for any $j\in \N$,
\begin{equation}\label{k_*}
    \|u_j \|_{L^{q_s}(\Omega)}^{q_s}\leq k_* (1+\|u_j \|_{\mathbb{H}_0^s(\Omega)}),
\end{equation}
for a convenient positive constant $k_*$.\\
Therefore, as a consequence of \eqref{phi k} and since $q<q_s$ we have
\begin{equation}
\begin{array}{ll}
k\geq \phi_{p,s}(u_j) & \geq\frac{1}{2}p_0 \|u_j\|_{ \mathbb{H}^s_0(\Omega)}^2 -\frac{\lambda }{2} \|u_j\|_{L^q(\Omega)}^q -\frac{1}{q_s} \|u_j\|_{L^{q_s}(\Omega)}^{q_s}\\
\\
& \geq \frac{1}{2}p_0 \|u_j\|_{ \mathbb{H}^s_0(\Omega)}^2 -C (1+\|u_j \|_{\mathbb{H}_0^s(\Omega)})^\frac{q}{q_s} -C (1+\|u_j \|_{\mathbb{H}_0^s(\Omega)}), 
\end{array}
\end{equation}
where $C$ is a positive constant.
Then combining this with \eqref{k_*}, we get that, for any $j\in \N$,
\begin{equation}
\|u_j\|_{ \mathbb{H}^s_0(\Omega)}^2 \leq \bar{k} \bigg[ (1+ \|u_j\|_{ \mathbb{H}^s_0(\Omega)} )+(1+\|u_j \|_{\mathbb{H}_0^s(\Omega)})^\frac{q}{q_s} \bigg],
\end{equation}
where $\bar{k}$ is a convenient positive constant. Since $q<q_s$, $(u_j)_{j\in\N}$ is bounded in $\mathbb{H}_0^s(\Omega)$.
\end{preuve} 
 \begin{lemme}
 Problem \eqref{formulation faible} admits a non trivial solution $u \in\mathbb{H}_0^s(\Omega). $
 \end{lemme}
 \begin{preuve}
     Since the sequence $(u_j)_{j\in\N}$ is bounded, and recalling that the fractional Sobolev space  $\mathbb{H}^s_0(\Omega)$ is a reflexive space, extract a subsequence, still denoted by $u_j$, so that 
 \begin{equation*}
 \begin{array}{ll}
 u_j \rightharpoonup u \hbox{ weakly in } \mathbb{H}^s_0(\Omega),
 \\
 u_j \rightharpoonup u \hbox{ weakly in } L^{q_s}(\Omega),
\\
 u_j \rightarrow u \hbox{ strongly in } L^q(\Omega) \hbox{ for any }q\in[1,q_s[,
\\
 u_j \rightarrow u \hbox{ a.e. on } \Omega.
 \end{array}
 \end{equation*}
 By taking into account the scalar product defined in \eqref{produit scalaire p}, for any $\varphi \in \mathbb{H}^s_0(\Omega)$ we have 
\begin{equation*}
     \displaystyle\int_{\R^{2n}} p(x)\frac{(u_j(x)-u_j(y))(\varphi(x)-\varphi(y))}{|x-y|^{n+2s}}dxdy 
     = \displaystyle\int_{\R^{2n}} p(x)\frac{(u(x)-u(y))(\varphi(x)-\varphi(y))}{|x-y|^{n+2s}}dxdy + o(1), 
\end{equation*}
as $j$ tends to $+\infty.$
We have also $|u_j|^{q_s-2} u_j \rightharpoonup |u|^{q_s-2}u$ weakly in $L^\frac{q_s}{q_s-1}(\Omega)$ as $j$ tends to $ +\infty,$
and 
$|u_j|^{q-2} u_j \rightarrow |u|^{q-2}u$ strongly in $L^\frac{q}{q-1}(\Omega)$ as $j$ tends to $+\infty.$\\
Moreover, for any $\varphi \in \mathbb{H}^s_0(\Omega),$ we have
\begin{equation}\label{lim phi'}
\begin{array}{ll}    

<\Phi'_{p,s}(u_j),\varphi> & = \displaystyle\int_{\mathbb{R}^n\times \mathbb{R}^n} p(x)\frac{(u_j(x)-u_j(y))(\varphi(x)-\varphi(y))}{|x-y|^{n+2s}}dxdy\\
\\
& -\lambda \displaystyle\int_{\Omega} |u_j(x)|^{q-2} u(x) \varphi(x) dx  - \displaystyle\int_{\Omega} |u_j(x)|^{q_s-2}u(x)\varphi(x) dx 
\end{array}
\end{equation}
By \eqref{conv 0},  $<\Phi'_{p,s}(u_j),\varphi> $ tends to zero as $j$ tends to $+\infty$. Thus, passing to the limit as $j$ tends to $+\infty$ in \eqref{lim phi'}, u satisfies
\begin{equation}
\begin{array}{ll}
\displaystyle\int_{\mathbb{R}^n\times \mathbb{R}^n} p(x)\frac{(u(x)-u(y))(\varphi(x)-\varphi(y))}{|x-y|^{n+2s}}dxdy\\
\\
 -\lambda \displaystyle\int_{\Omega} |u(x)|^{q-2} u(x) \varphi(x) dx  - \displaystyle\int_{\Omega} |u(x)|^{q_s-2}u(x)\varphi(x) dx = 0,
\end{array}
\end{equation}
for any $\varphi \in \mathbb{H}^s_0(\Omega)$; that is, u is a solution of problem \ref{formulation faible}, and the second point follows.\\
Now let us prove that the solution $u$ is not zero.
Let us suppose by contradiction that $u \equiv 0$ in $\Omega$.\\
Since $u \in \mathbb{H}^s_0(\Omega)$ $u \equiv 0$ in $\mathbb{R}^n$ . We know that $(u_j)_{j\in\N}$ is bounded in $ \mathbb{H}^s_0(\Omega)$, and
\begin{equation}
\begin{array}{ll}    
 <\Phi'_{p,s}(u_j),u_j> & = \displaystyle\int_{\mathbb{R}^n\times \mathbb{R}^n} p(x)\frac{|u_j(x)-u_j(y)|^2}{|x-y|^{n+2s}}dxdy\\
 \\
& -\lambda \displaystyle\int_{\Omega} |u_j(x)|^{q} dx  - \displaystyle\int_{\Omega} |u_j(x)|^{q_s} dx,
\end{array}
\end{equation}
Since \eqref{conv 0} holds true for any $\varphi \in \mathbb{H}^s_0(\Omega)$,  $<\Phi'_{p,s}(u_j),u_j>$ tends to zero as $j$ tends to $+\infty$.
Using also $\displaystyle\int_{\Omega} |u_j(x)|^{q} dx \rightarrow 0 $ as $j$ tends to $+\infty$, we write
\begin{equation}\label{diff tend 0}
\displaystyle\int_{\mathbb{R}^n\times \mathbb{R}^n} p(x)\frac{|u_j(x)-u_j(y)|^2}{|x-y|^{n+2s}}dxdy- \displaystyle\int_{\Omega} |u_j(x)|^{q_s} dx \rightarrow 0, \hbox{ as } j \rightarrow +\infty.
\end{equation}
We know that the sequence $(\|u_j \|_{\mathbb{H}^s_0(\Omega)})_{j\in\N}$ is bounded in $\R$. \\
Hence, the sequence $\bigg(\displaystyle\int_{\mathbb{R}^n\times \mathbb{R}^n} p(x)\frac{|u_j(x)-u_j(y)|^2}{|x-y|^{n+2s}}dxdy\bigg)_{j\in\N}$ is also bounded in $\R$ since the weight $p$ is bounded. Extract a subsequence, we assume that
\begin{equation}\label{p tend l}
\displaystyle\int_{\mathbb{R}^n\times \mathbb{R}^n} p(x)\frac{|u_j(x)-u_j(y)|^2}{|x-y|^{n+2s}}dxdy \rightarrow l,
\end{equation}
as $j$ tends to $+\infty,$ and by \eqref{diff tend 0} we easily deduce that
\begin{equation}\label{u_j tend l}
\displaystyle\int_{\Omega} |u_j(x)|^{q_s} dx \rightarrow l,
\end{equation}
as $j$ tends to $+\infty.$\\
Moreover, by \eqref{conv c}, we have
\begin{equation}
    \begin{array}{ll}
 \frac{1}{2} \displaystyle\int_{\mathbb{R}^n \times \mathbb{R}^n} p(x) \frac{|u_j(x)-u_j(y)|^2}{|x-y|^{n+2s}}dx dy   - \frac{\lambda}{q} 
 \displaystyle\int_{\Omega}|u_j(x)|^q dx -\frac{1}{q_s} \displaystyle\int_{\Omega}|u_j(x)|^{q_s}dx \rightarrow c,
    \end{array}
\end{equation}
as $j$ tends to $+\infty$. Since  $\displaystyle\int_{\Omega} |u_j(x)|^{q} dx \rightarrow 0 $ as $ j $ tends to $+\infty$, using \eqref{p tend l} and \eqref{u_j tend l}, we deduce that
\begin{equation}\label{égalité c}
c=(\frac{1}{2}-\frac{1}{q_s})l=\frac{s}{n}l
\end{equation}
Using Proposition \ref{vérif 2.11} $c\geq \beta >0$, it is obvious that $l>0$. Furthermore, by definition of $S_s$, we have
\begin{equation}
\displaystyle\int_{\mathbb{R}^n\times \mathbb{R}^n} p(x)\frac{|u_j(x)-u_j(y)|^2}{|x-y|^{n+2s}}dxdy \geq p_0 S_s \| u_j \|^2_{L^{q_s}(\Omega)}.
\end{equation}
Passing to the limit as $ j $ tends to $+\infty$ and combining with \eqref{p tend l} and \eqref{u_j tend l}, we get
\begin{equation}
l\geq p_0 S_s l^{\frac{2}{q_s}},
\end{equation}
and taking into account \eqref{égalité c}, we have
\begin{equation}
    c\geq\frac{s}{n} p_0^{\frac{n}{2s}}  S_s^{\frac{n}{2s}},
    \end{equation}
a contradiction to the fact that $c <\frac{s}{n} p_0^{\frac{n}{2s}} S_s^{\frac{n}{2s}}$. Thus $u$ is not trivial and the proof of Theorem \ref{prob sans minimisation} is complete.\\
\end{preuve}





\begin{thebibliography}{Zz99}
 
\bibitem[1] {AR} A.  Ambrosetti,  P. Rabinowitz, Dual variational methods in critical point theory and applications, 
J. Funct. Anal. \textbf{14}, (1973), 349-381.

\bibitem[2] {A} T. Aubin, Equations différentielles non linéaires et problème de Yamabe concernant la courbure scalaire, 
J. Math. Pures Appl. (9), no. 3, \textbf{55} (1976), 269-296.

 \bibitem[3]{AB} B. Abdellaoui, R. Bentifour,  Caffarelli-Kohn-Nirenberg type inequalities of fractional order with applications,Journal of Functional Analysis, 10, \textbf{272} (2017), 3998-4029. 

\bibitem[4] {ABS} L. Appolloni, G.M. Bisci, S. Secchi, On critical Kirchhoff problems driven by the fractional Laplacian, 
Calculus of Variations and Partial Differential Equations, 6, \textbf{60} (2021), 209.

\bibitem[5]{BHY} S. Bae, R. Hadiji, H. Yazidi, Nonlinear existence result for a quasi-linear elliptic PDE, J. Math. Anal. Appl 396.1 (2012), 98-107.


\bibitem[6]{BH} A. Beaulieu, R. Hadiji, Remarks on solutions of a fourth-order problem, 
Applied Mathematics Letters, {\bf 19} (2006), 661-666.

\bibitem[7]{BHY} A. Benhamida, R. Hadiji, H. Yazidi, Existence results for elliptic equation  involving polyharmonic operator and a  critical  growth, to appear.



\bibitem[8] {BRS} G.M. Bisci, V.D. R{\u{a}}dulescu, R. Servadei, Variational methods for non local fractional problems, Cambridge University Press, \textbf{162} (2016).

 \bibitem[9]{B} H. Brezis, Some variational problems with lack of compactness,
 Proc. Symp. Pure Math. Vol, {\bf 45} Part 1 (F. Browder ed.), Amer. Math. Soc (1986), 165-201.

\bibitem[10]{BL} H. Brezis, E. Lieb, A relation between pointwise convergence of functions and convergence of functionals, 
Proc. Amer. Math. Soc, {\bf 88} (1983), 486-490.

\bibitem[11]{BN1} H. Brezis, L. Nirenberg, Positive solutions of nonlinear elliptic equations involving critical Sobolev exponents, 
Commun. Pure Appl. Math, {\bf 36} (1983), 437-477.


\bibitem[12]{CG} S.-Y. A. Chang and M. d. M. González, Fractional Laplacian in conformal geometry, Adv. Math.
226:2 (2011), 1410–1432.

\bibitem[13]{CS} L. Caffarelli, L. Silvestre,  An extension problem related to the fractional Laplacian, 
Communications in partial differential equations, 8, \textbf{32} (2007), 1245-1260.

\bibitem[14]{CT} A. Cotsiolis, N.K. Tavoularis, Best constants for Sobolev inequalities for higher order fractional derivatives, 
Journal of mathematical analysis and applications, 1, \textbf{295} (2004), 225-236. 

\bibitem[15]{Es} J. F. Escobar, “Conformal deformation of a Riemannian metric to a scalar flat metric with constant mean curvature on the boundary”, Ann. of Math. (2) 136:1 (1992), 1-50. 

\bibitem[16]{FLS} R. Frank, E. Lieb, R. Seiringer, Hardy-Lieb-Thirring inequalities for fractional Schr{\"o}dinger operators, 
Journal of the American Mathematical Society, 4, \textbf{21} (2008), 925-950.

\bibitem[17]{FV} F. Ferrari, I. Verbitsky, Radial fractional Laplace operators and Hessian inequalities, 
Journal of Differential Equations, 1, \textbf{253} (2012), 244-272.

\bibitem[18]{GQ} M. Gonz{\'a}lez Nogueras, J. Qing, Fractional conformal Laplacians and fractional Yamabe problems.,
Analysis \& PDE, 7, \textbf{6} (2013), 1535-1576.

\bibitem[19]{H} R. Hadiji,  A nonlinear problem with a weight and a non-vanishing boundary datum, 
Pure and Applied Functional Analysis, 4, \textbf{5} (2020), 965-980.


\bibitem[20]{HMPY} R. Hadiji, R. Molle, D. Passaseo, H. Yazidi, Localization of solutions for nonlinear elliptic problems with critical growth, Comptes rendus. Mathématique 343.11-12 (2006): 725-730.

\bibitem[21] {HV} R. Hadiji, F. Vigneron, Existence of solutions of a nonlinear eigenvalue problem with a variable weight, 
Journal of Differential Equations, 2-3, \textbf{266} (2019), 1488-1513.

\bibitem[22]{HY} R. Hadiji and  H. Yazidi, Problem with critical Sobolev exponent and with weight, 
Chinese Annal. Math. ser. B, 3, \textbf{28}, (2007), 327-352.

\bibitem[23]{L} E. Lieb, Sharp constants in the Hardy-Littlewood-Sobolev and related inequalities, 
 Ann. of Math. 2, \textbf{118} (1983), 349-374.

\bibitem[24]{LP} J.M. Lee, T.H. Parker, The Yamabe problem, 
Bull. Amer. Math. Soc., no. 1, \textbf{17} (1987), 37-91.

\bibitem[25]{MOSW} W. Magnus, F. Oberhettinger, R.P. Soni, E. Wigner. Formulas and theorems for the special functions of mathematical physics, Physics Today, 12, \textbf{20} (1967), 81-83. 

\bibitem[26]{NPV} E. Di Nezza, G. Palatucci, E. Valdinoci, Hitchhiker's guide to the fractional Sobolev spaces, 
Bulletin des sciences mathématiques, 5, \textbf{136} (2012), 521-573.

\bibitem[27]{NS} H. Nguyen, M. Squassina, Fractional Caffarelli-Kohn-Nirenberg inequalities, 
Journal of Functional Analysis, 9, \textbf{274} (2018), 2661-2672.

\bibitem[28]{S} R. Servadei. The Yamabe equation in a nonlocal setting, 
Advances in Nonlinear Analysis, 3, \textbf{3} (2013), 235-270.

\bibitem[29]{SV} R. Servadei, E. Valdinoci, Mountain pass solutions for nonlocal elliptic operators, Journal of Mathematical Analysis and Applications, 2, \textbf{389} (2012), 887-898.











 
 \end{thebibliography}
\end{document}